\documentclass{amsart}
\usepackage{amsfonts,amssymb,amsmath,a4wide}

\newtheorem{teo}{Theorem}[section]
\newtheorem{lema}{Lemma}[section]
\newtheorem{pro}{Proposition}[section]
\newtheorem{defi}{Definition}[section]
\newtheorem{rema}{Remark}[section]

\numberwithin{equation}{section}

\def \z{\zeta}

\def \g{\mathfrak{g}}
\def \K{\mathcal{K}}
\def \hook{\lrcorner}
\def\sideremark#1{\ifvmode\leavevmode\fi\vadjust{\vbox to0pt{\vss
 \hbox to 0pt{\hskip\hsize\hskip1em
 \vbox{\hsize2.5cm\tiny\raggedright\pretolerance10000
 \noindent #1\hfill}\hss}\vbox to8pt{\vfil}\vss}}}%

\begin{document}\larger[2]
\title{Skew-symmetric prolongations of Lie algebras and applications}
\author[P.-A.Nagy]{Paul-Andi Nagy}
\address[P.-A. Nagy]{Department of Mathematics, University of Auckland, Private Bag 92019, Auckland, New Zealand}
\email{nagy@math.auckland.ac.nz}
\date{\today}
\begin{abstract}We study the skew-symmetric prolongation of a Lie subalgebra $\g \subseteq \mathfrak{so}(n)$, in other words the intersection 
$\Lambda^3 \cap (\Lambda^1 \otimes \g)$.
We compute this space in full generality. Applications include uniqueness results for connections with skew-symmetric torsion and also the proof of the Euclidean version of a 
conjecture posed in \cite{ofarill} concerning
a class of Pl\"ucker-type embeddings. We also derive a classification of the metric k-Lie algebras (or Filipov algebras), in positive signature and 
finite dimension. Prolongations of Lie algebras can also be used to finish the classification, started in \cite{datri}, of manifolds admitting Killing frames, or equivalently flat connections 
with $3$-form torsion. Next we study specific properties of invariant $4$-forms of a given metric representation and apply these considerations to classify 
the holonomy representation of metric connections with vectorial torsion, that is with torsion contained in $\Lambda^1 \subseteq \Lambda^1 \otimes \Lambda^2$. 
\end{abstract}
\maketitle
\tableofcontents
\section{Introduction}
Let $(M^n,g)$ be an oriented Riemannian manifold of dimension $n$. A $G$-structure on $M$ consists in a reduction of the frame bundle of $(M^n,g)$ to a Lie subgroup 
$G \subseteq SO(n)$. While this is an interesting 
object by itself, it appears naturally in the context of many partial differential equations where it is often additionally equipped with a connection 
preserving the $G$-structure, generally with non-vanishing torsion. Without further assumptions this is known to be fairly unrestrictive \cite{Hano}. Sometimes the 
reduction of the frame bundle of $(M^n,g)$ is realised by means of the isotropy group of a tensor field on $M$ of specific algebraic type, as required for instance 
by models in string theory \cite{stro}. When combining the classification of $G$-structures with the study of the torsion tensor of a given $G$-connection
(see \cite{Cartan,GHer,cs} for instance)  it is possible to distinguish several types of geometries. An approach to the classification of each type could be 
through the study of the holonomy of the $G$-connection involved. 

In the case of connections of Levi-Civita or torsion-free type powerful classification results  
(see \cite{br1,br2,MS}) are available and each of the resulting geometries e.g. K\"ahler, Calabi-Yau, quaternion-K\"ahler, Joyce manifolds is in the mainstream 
of current research. The theory of connections with non-vanishing torsion is less well established 
and many classical results from Riemannian holonomy theory fail to hold, such as the deRham splitting theorem.

In this paper we shall mainly consider metric connections with skew-symmetric or vectorial torsion on $(M^n,g)$, meaning that their torsion is 
assumed to be in $\Lambda^3M$ or $\Lambda^1M \subseteq \Lambda^1M
\otimes \Lambda^2M$, respectively. These classes are labeled 
$W_1$ and $W_3$ and shall be given a more detailed description later on. There are situations when the class $W_1$ emerges naturally from a specific class of 
$G$-manifolds such as nearly-K\"ahler \cite{Gray,nagy} or $Spin(7)$-structures  \cite{Iv}, for instance. 
Concerning the class $W_1$ we wish to address the question of uniqueness when the holonomy representation is fixed, that is if the requirement on a connection $D$ in the 
class $W_1$ to have holonomy contained in a fixed Lie subgroup 
$G \subset SO(n)$ implies its uniqueness. Note that such a connection might not always exist, and that the uniqueness problem has been solved for those representations 
arising as Riemannian, non-locally symmetric holonomies in 
\cite{Iv,Fri1, Fri2}, as well as for the irreducible representations of $SO(3)$ and $PSU(3)$ in dimensions $5$ and $8$ respectively  \cite{BoNu, nurow}. To treat the general case we will consider, for a 
given faithful representation $(\g,V)$ of a Lie algebra on an Euclidean vector its skew-symmetric prolongation 
\begin{equation*} \label{pdef}
\Lambda^3V \cap (\Lambda^1V \otimes \g)=\{T \in \Lambda^3V : X \lrcorner T \in \g, \ \mbox{for all} \ X \in V\}.
\end{equation*}
We are able to compute explicitly this space as follows.
\begin{teo} \label{res1}
Let $(\g,V)$ be a faithful and irreducible representation of a Lie algebra onto an Euclidean space. If $\g$ is proper then either:
\begin{itemize}
\item[(i)]$\Lambda^3V \cap (\Lambda^1V \otimes \g)=\{0\}$;
\item[] or
\item[(ii)]
$(\g,V)$ is the adjoint representation of a compact, simple Lie algebra, when $\Lambda^3V \cap (\Lambda^1V \otimes \g)$ is $1$-dimensional, generated by the 
Cartan form of $\g$. 
\end{itemize}
\end{teo}
We also show how this can be used to compute the skew-symmetric prolongation of an arbitrary orthogonal representation.
As a direct application of Theorem \ref{res1} we obtain that:
\begin{teo} \label{uniq1}
Let $(M^n,g)$ be a connected Riemannian manifold equipped with a Riemannian $G$-structure. Denote by $(G,V)$ the corresponding representation on tangent spaces and assume it is irreducible. Then:
\begin{itemize}
\item[(i)] if a metric connection $D$ with 
$3$-form torsion which respects the $G$ structure exists it must be unique, exception made of the case when  $(G,V)$ is the adjoint representation of 
a compact, simple Lie group; 
\item[(ii)] the latter case occurs if and only if $(M^n,g)$ carries a $3$-form $t$ in $\Lambda^3M$ such that $Dt=0$ and $t^2=1$ in the Clifford algebra bundle of $M$.
\end{itemize}
\end{teo}
Theorem \ref{res1} has a few other applications we shall briefly describe now. The first consists in the explicit description of the holonomy of connections with constant 
torsion in flat space, or equivalently of holonomy algebras generated by a $3$-form, a study initiated in \cite{agfr}. We prove:
\begin{teo} \label{hol3} Let $(V,g)$ be an Euclidean vector space and let $T\neq 0$ belong to $\Lambda^3V$. If the subalgebra $\g^{\star}_T=Lie \{(X \lrcorner T)^{\sharp} : X \in V\}$ of $\mathfrak{so}(V)$ acts irreducibly on 
$V$ then either
\begin{itemize}
\item[(i)] $\g_T^{\star}=\mathfrak{so}(V)$
\item[] or 
\item[(ii)] $(\g_T^{\star},V)$ is the adjoint representation of a compact simple Lie algebra with $T$ is proportional to its Cartan form.
\end{itemize}
\end{teo}
Here ${}^{\sharp} : \Lambda^2V \to \mathfrak{so}(V)$ is the linear isomorphism induced by the metric $g$ above. Closely related is the problem (see \cite{datri}) of classifying Riemannian manifolds $(M^n,g)$ which admit a Killing frame, that is 
an orthonormal frame $\{\z_i, 1 \le i \le n\}$ such that $\z_i$ is a Killing vector field  for the metric $g$ for all $1 \le i \le n$. The existence of a Killing frame implies that $(M^n,g)$ carries a flat metric connection with $3$-form torsion and the converse 
is true if $M$ is simply connected (see \cite{datri}). 
Examples include flat spaces, Lie groups equipped with a bi-invariant metric and also $S^7$ with its canonical metric. In the latter case a Killing frame is obtained by restricting the octonian multiplication to $S^7$. Moreover, it was 
shown in \cite{datri} that apart from these examples the only simply connected, complete Riemannian manifold which might carry a Killing frame is the symmetric space $SU(2m) \slash Sp(m), m \ge 2$. Making use of Theorem \ref{hol3} we are 
able to exclude this instance and prove:
\begin{teo} \label{class3flat}
Let $(M^n,g)$ be simply connected and complete. If it admits a flat metric connection with $3$-form torsion or, equivalently, a Killing frame then $(M^n,g)$ is a Riemannian product of with factors in one of the following classes:
\begin{itemize}
\item[(i)] Euclidean spaces;
\item[(ii)] simply connected, compact semisimple Lie groups equipped with a bi-invariant metric;
\item[(iii)]  $S^7$ with a metric of constant sectional curvature.
\end{itemize}
\end{teo}
Theorem \ref{res1} can also be used to study a class of Pl\"ucker type relations \cite{ofarill}, having their origins in the study of maximally supersymmetric solutions 
of ten and eleven dimensional supergravity theories \cite{ofarpapa}. Note that these require the flatness of a spinorial connection obtained by modifying the Levi-Civita connection by a differential form of arbitrary degree. 
We confirm the Euclidean case of the conjecture in \cite{ofarill} by proving:
\begin{teo} \label{plckclass} Let $(V^n,g)$ be a Euclidean vector space and let $T \neq 0$ in $\Lambda^pV$ satisfy  
$$ [\z_1 \lrcorner \ldots \z_{p-2} \lrcorner T,T]=0
$$
for all $\z_i, 1 \le i \le p-2$ in $V$. If $p \ge 4$ the form $T$ is decomposable i.e it is an orthogonal sum of simple forms.
\end{teo}
Here the reader is referred to Section 2 for unexplained notation. We also obtain a classification of the so-called $n$-Lie (or Filipov) algebras, starting from the observation \cite{ofarill} that those are in 
$1:1$ correspondence with the Pl\"ucker relations mentioned above. In the case of neutral signature, examples where Theorem \ref{plckclass} fails to hold for a class of $4$-forms have been constructed in \cite{Cortes}. 

In the second part of the paper we classify the holonomy of metric connections with vectorial torsion. More precisely we look at Riemannian manifolds $(M^n,g)$ equipped 
with a metric connection $D$ of the form 
$$ D_X=\nabla_X+X \wedge \theta
$$
for all $X$ in $TM$, where $\nabla$ is the Levi-Civita connection of the metric $g$ and $\theta$ is a $1$-form on $M$. The triple $(M^n,g,D)$ will be called {\it{closed}} 
if $\theta$ is closed, that is $d\theta=0$. The main class of examples in this direction is obtained when $(M^n,g)$ has a locally parallel structure (l.c.p. for short).
That is we require $g$ not to be locally conformally symmetric and moreover it must be locally conformal to a metric admitting a Riemannian holonomy 
reduction (see \cite{redbk} for details). Note that such a metric needs not be globally conformally parallel.
For results and examples in this direction we refer the reader to \cite{IvPar,Ornea,OrneaPicini} and references therein. We prove:

\begin{teo} \label{res2} Let $(M^n,g)$ be a connected and oriented Riemannian manifold and let $D$ be a metric connection with vectorial torsion. If the 
holonomy representation $(G,V)$ of $D$ is irreducible with $G$ proper in $SO(n)$ then:
\begin{itemize}
\item[(i)] $(M^n,g,D)$ is closed. Moreover, either:
\item[(ii)] $g$ is conformal to a non-flat locally symmetric Riemannian metric on $M$,
\item[(iii)] the connection $D$ is flat,
\item[or]
\item[(iv)] $g$ admits a l.c.p. structure.
\end{itemize}
\end{teo}
Note that closedeness for $(M^n,g,D)$ as above equally holds \cite{agfrivec} in case $D$ admits a non-zero parallel spinor.
The paper is organised as follows. In Section 2 we start by presenting some well known facts related to the classification of connections with torsion and also to 
some of their basic properties including the Bianchi identity and the symmetries of the non-Riemannian curvature tensors. Next, and by following 
\cite{cs} essentially,  we recall the results on Berger algebras and algebraic curvature tensors we will need in what follows.  

Section 3 is devoted to the proof of Theorem \ref{res1} which is build on the observation that the second symmetric tensor power of the skew-symmetric 
prolongation space maps into the space of algebraic curvature tensors. This is initially done for irreducible representations and we also apply the results 
to prove Theorem \ref{hol3}. Finally we extend Theorem \ref{res1} to the case of an arbitrary 
representation, by using a decomposition algorithm from \cite{cs}. In Section 4 we present the proof of Theorem \ref{uniq1} whilst in Sections 5 and 6 we prove Theorems \ref{class3flat}, \ref{plckclass} again 
by reduction to a Lie algebra prolongation problem. Section 7 contains the necessary ingredients for the proof of Theorem \ref{res2}. We use the construction of invariant $4$-forms for an orthogonal metric representation 
given in \cite{kos} (see also \cite{alco}) and isolate the algebraic properties which are needed in order to prove closedeness. This enables us to use results from \cite{agfrivec} and \cite{cs} in order to obtain the  classification. 
Finally, in Section 8 we look at a class of representations which occur naturally when 
one considers the existence problem for metric connections with skew-torsion. Those representations were classified in \cite{Fri2} under the requirement that their Lie algebra has the involution 
property. Through arguments similar to those used in the proof of 
Theorem \ref{res1}, we show that the involutivity assumption can be removed. 
\section{Metric connections with torsion}
Let $(M^n,g)$ be an oriented Riemannian manifold. To the metric $g$ one attaches its Levi-Civita connection $\nabla$ which leaves $g$ invariant and moreover it is 
torsion free. More generally, 
a metric connection $D$ on $(M^n,g)$ is a linear connection on $TM$ preserving the Riemannian metric, that is 
$Dg=0$. Any such connection $D$ can be written as $D=\nabla+\eta$ where the tensor $\eta$ belongs to 
$\Lambda^1M \otimes \Lambda^2M$. The torsion tensor $Tor$ in $\Lambda^2M \otimes \Lambda^1M$ of the connection $D$ is given by 
\begin{equation*}
\begin{split}
Tor(X,Y)=&D_XY-D_YX-[X,Y]\\
=&\eta_XY-\eta_YX
\end{split}
\end{equation*}
for all $X,Y$ in $TM$. Moreover, the connection $D$ is extended to $\Lambda^{\star}M$ by defining 
\begin{equation*}
D_X\varphi=\nabla_X \varphi+[\eta_X, \varphi]
\end{equation*}
for all $X$ in $TM$ and whenever $\varphi$ belongs to $TM$. Here we have defined the commutator 
\begin{equation*}
[\alpha, \gamma]=\sum \limits_{i=1}^n (e_i \lrcorner \alpha) \wedge (e_i \lrcorner \gamma)
\end{equation*}
for all $\alpha$ in $\Lambda^2M$ and for all $\gamma$ in $\Lambda^{\star}M$, where $\{e_i, 1 \le i \le n\}$ is some 
local orthonormal basis in $TM$. In what follows we shall identify $\Lambda^2M$ and $\mathfrak{so}(TM)$, by writing any $2$-form $\alpha$ as $\alpha=g(F \cdot, \cdot)$ for 
some skew-symmetric endomorphism of $TM$. Also, we shall systematically use the metric $g$ to identify $\Lambda^1M$ to $TM$. For further computations, it is also useful 
to note that in some local
orthonormal basis $\{ e_i, 1 \le i \le n\}$ of $TM$, one has 
$$g([F,G] \cdot, \cdot)=\sum \limits_{i=1}^n (e_i \lrcorner \alpha) \wedge (e_i \lrcorner \beta)$$
for all $F,G$ in $\mathfrak{so}(TM)$ with dual  $2$-forms $\alpha=g(F \cdot, \cdot)$ and $\beta=g(G \cdot, \cdot)$, where the bracket in l.h.s. is the usual 
Lie bracket of $\mathfrak{so}(TM)$.

Consider now the Cartan decomposition 
\begin{equation*}
\Lambda^1M \otimes \Lambda^2M=W_1 \oplus W_2 \oplus W_3 
\end{equation*}
into irreducible components under the action of $SO(n)$. Explicitly, one has 
\begin{equation*}
\begin{split}
W_1=&\Lambda^3M\\
W_2=& (\Lambda^1M \otimes \Lambda^2M) \cap Ker a \cap Ker t\\
W_3=&\Lambda^1M.
\end{split}
\end{equation*}
Here, for any $p \ge 0$ the total alternation map $a : \Lambda^1M \otimes \Lambda^pM \to \Lambda^{p+1}M$ and the trace map 
$t : \Lambda^1M \otimes \Lambda^pM \to \Lambda^{p-1}M$ are given, respectively, by 
\begin{equation*}
\begin{split}
a(\eta)=& \sum \limits_{i=1}^n e_i \wedge \eta_{e_i} \\
t(\eta)=&  \sum \limits_{i=1}^n e_i \lrcorner \eta_{e_i}
\end{split}
\end{equation*}
for all $\eta$ in $\Lambda^1M \otimes \Lambda^pM$, where $\{e_i, 1 \le i \le n\}$ is some local orthonormal frame on $M$. We also note that the embedding of 
$\Lambda^1M$ into $\Lambda^1M \otimes \Lambda^2M$ is given by $X \mapsto X \wedge \cdot$.
Since the torsion tensor of any metric connection on $(M^n,g)$ lives in $\Lambda^2M \otimes \Lambda^1M$ it follows \cite{Cartan} that there are nine 
main classes to be considered. In this paper we shall be mainly interested in the class $W_1$, when the torsion is given by a $3$-form on $M$ and in the class 
$W_3$ when the torsion tensor is determined by a $1$-form. 

Before giving details, we recall that any metric connection $D$ has its curvature tensor $R^D$ in $\Lambda^2M \otimes \Lambda^2M$ defined by $R^D(X,Y)=-D^2_{X,Y}+D^2_{Y,X}-D_{Tor(X,Y)}$ for all 
$X,Y$ in $TM$. In the case of the Levi-Civita connection it will be simply denoted by $R$. When $D=\nabla+\eta$ for some $\eta$ in $\Lambda^1M \otimes \Lambda^2M$ 
one obtains the following comparison formula 
\begin{equation} \label{Kcomp}
R^D(X,Y)=R(X,Y)-\biggl [(D_X\eta)_Y-(D_Y\eta)_X \biggr ]+[\eta_X, \eta_Y]-\eta_{Tor(X,Y)}
\end{equation}
for all $X,Y$ in $TM$. Let us also recall that the Bianchi operator $b_1 : \Lambda^2M \otimes \Lambda^2M \to \Lambda^1M \otimes \Lambda^3M$ is defined by 
\begin{equation*}
(b_1Q)_X=\sum \limits_{i=1}^n e_i \wedge Q(e_i, X)
\end{equation*}
for all $(X,Q)$ in $TM \times (\Lambda^2M \otimes \Lambda^2M)$, where  $\{e_i, 1 \le i \le n\}$ is some local orthonormal frame on $M$. 
\subsection{$3$-form torsion}
Explicitly, a metric connection $D$ belongs to the class $W_1$ if and only if it is of the form 
\begin{equation*}
D_X=\nabla_X+\frac{1}{2}T_X
\end{equation*}
for all $X$ in $TM$, where $T$ belongs to $\Lambda^3M$. Here and henceforth, for any $T$ in $\Lambda^3M$ and for any $X$ in $TM$ we denote by $T_X$ in $\mathfrak{so}(TM)$ the dual 
of $X \lrcorner T$ w.r.t. to the metric $g$, that is $X \lrcorner T=g(T_X \cdot, \cdot)$. In this case one has $Tor=T$. Specific objects of relevance are the $4$-form $\Omega^T$ in $\Lambda^4M$ given by 
\begin{equation*}
\Omega^T=\frac{1}{2} \sum \limits_{i=1}^n T_{e_i} \wedge T_{e_i}
\end{equation*}
for some arbitrary local orthonormal frame $\{e_i, 1 \le i \le n\}$ on $M$. It satisfies 
\begin{equation} \label{locomega}
\Omega^T(X,Y)=-[T_X,T_Y]+T_{T_XY}
\end{equation}
where the notation $\Omega^T(X,Y)=Y \lrcorner X \lrcorner \Omega^T$ for all $X,Y$ in $TM$ shall be used constantly in the subsequent. Secondly, we define the tensor $R^T$ in $\Lambda^2M \otimes \Lambda^2M$ by 
\begin{equation*}
R^T(X,Y)=[T_X,T_Y]+2T_{T_XY}
\end{equation*}
for all $X,Y$ in $TM$. At this point we recall that the Ricci contraction operator $Ric : \otimes^2\Lambda^2M \to \otimes^2M$ is given by 
\begin{equation*}
Ric(Q)(X,Y)=\sum \limits_{i=1}^{n}Q(X,e_i,Y,e_i)
\end{equation*}
whenever $Q$ belongs to $\otimes^2\Lambda^2M$ and for all $X,Y$ in $TM$, where $\{e_i, 1 \le i \le n\}$ is some local orthonormal frame on $M$. Note that the Ricci operator preserves 
tensor type w.r.t the splitting $\Lambda^2M \otimes \Lambda^2M=S^2(\Lambda^2M) \oplus \Lambda^2(\Lambda^2M)$. The following is proved by a routine verification which is left to the reader. 
\begin{lema} \label{curvT1} The following hold:
\begin{itemize}
\item[(i)] $R^T$ is an algebraic curvature tensor, that is $b_1R^T=0$;
\item[(ii)] $Ric(R^T)(X,Y)=3\langle X \lrcorner T, Y \lrcorner T \rangle $ for all $X,Y$ in $TM$.
\end{itemize}
\end{lema}
Together with $D$ comes its associated exterior derivative $d_D : \Lambda^{\star}M \to \Lambda^{\star+1}M$ given by 
\begin{equation*}
d_D=\sum \limits_{i=1}^n e_i \wedge D_{e_i}
\end{equation*}
where $\{e_i, 1 \le i \le n\}$ is some local orthonormal frame on $M$. It is related to the ordinary exterior derivative $d$ by 
\begin{equation*}
d_D\varphi=d\varphi-T \bullet \varphi
\end{equation*}
where we have defined 
$$ T \bullet \varphi = \sum_i (e_i \hook T) \wedge (e_i \hook \varphi), $$
whenever $\varphi$ belongs to $\Lambda^{\star}M$ and $\{e_i, 1 \le i \le n\}$ is some local orthonormal frame on $M$. Note that $\Omega^T=2T \bullet T$. 

Finally, let $\Theta$ in $\Lambda^2(\Lambda^2M)$ be given by 
\begin{equation*}
\Theta(X,Y)=(D_XT)_Y-(D_YT)_X-\frac{1}{2}d_DT(X,Y)
\end{equation*} 
for all $X,Y$ in $TM$. The comparison formula \eqref{Kcomp} becomes now 
\begin{equation} \label{curv3form}
R^D=(R-\frac{1}{12}R^T)-\frac{1}{2}\Theta-(\frac{1}{4}d_DT+\frac{1}{3}\Omega^T).
\end{equation}
It identifies the components of $R^D$ along the orthogonal splitting $\otimes^2\Lambda^2M=\mathcal{K}M \oplus \Lambda^2(\Lambda^2M) \oplus \Lambda^4M$ where the bundle of algebraic curvature tensors on $M$ is defined by 
$ \mathcal{K}M=\otimes ^2\Lambda^2M \cap Kerb_1$. A few immediately useful consequences of \eqref{curv3form} are listed above.
\begin{lema} \label{bianchis}
The following hold whenever $X$ belongs to $TM$:
\begin{itemize}
\item[(i)]  $(b_1\Theta)_X=-2D_XT+\frac{1}{2}X \lrcorner d_DT$;
\item[(ii)] $(b_1R^D)_X=D_XT+\frac{1}{2}X \lrcorner dT$.
\end{itemize}
\end{lema}
\begin{proof} (i) follows by direct computation whilst (ii) is an immediate consequence of (i) and \eqref{curv3form}.
\end{proof}
For future use it is also useful to consider the $1$-parameter family of metric connections 
$$ D^t_X=\nabla_X+\frac{t}{2}T_X
$$ 
for all $X$ in $TM$, where $t$ belongs to $\mathbb{R}$. Then we recover the Levi-Civita connection as well as $D$ from $D^0=\nabla$ and $D^1=D$, respectively. Moreover 
the torsion type is preserved, since $Tor(D^t)=tT$ for all $t$ in $\mathbb{R}$. This type of deformation has been used in various contexts in \cite{bismut, kostant2, agfr}.
We denote by $R^t$ the curvature tensor of the connection $D^t, t \in \mathbb{R}$ and 
we shall work towards obtaining a comparison formula of those tensors
with the reference curvature tensor $R^D$.
\begin{pro} \label{curvpar}
For any $t$ in $\mathbb{R}$ we have: 
\begin{equation} \label{param} 
R^ t=R^D-\frac{t-1}{2} \Theta-\frac{t-1}{4}d_DT-\frac{t^2-1}{12}R^T+\frac{(t-1)(t-2)}{6}\Omega^T
\end{equation}
\end{pro}
\begin{proof}
From \eqref{curv3form}, applied for the connection $D^t, t$ in $\mathbb{R}$ we get 
\begin{equation*} \label{Kcompt}
R^t(X,Y)=R(X,Y)-\frac{t}{2} \biggl [ (D^t_XT)_Y-(D^t_YT)_X\biggr ]-t^2(\frac{1}{12}R^T(X,Y)+\frac{1}{3}\Omega^T(X,Y))
\end{equation*}
for all $X,Y$ in $TM$. Now
\begin{equation*}
D^t_XT=D_XT+\frac{t-1}{2}[T_X,T]=D_XT-\frac{t-1}{2}X \lrcorner \Omega^T
\end{equation*}
after observing that  $X \lrcorner \Omega^T = - [T_X, T]$ for all $X$ in $TM$. By also using \eqref{curv3form} this time for $D$ it follows after comparison that  
\begin{equation*}
\begin{split}
R^t(X,Y)=&R^D(X,Y)-\frac{t-1}{2} \biggl [ (D_XT)_Y-(D_YT)_X \biggr ]+\frac{t(t-1)}{2} \Omega^T(X,Y)\\
&-(t^2-1)(\frac{1}{12}R^T(X,Y)+\frac{1}{3}\Omega^T(X,Y))
\end{split}
\end{equation*}
whenever $X,Y$ belong to $TM$.The proof of the claim follows now immediately by using the definition of $\Theta$. 
\end{proof}
\subsection{Vectorial torsion}
We shall summarise now some of the facts and results we will need concerning connections in the class $W_3$. If the connection $D$ is in the class $W_3$ it can be written 
as 
\begin{equation*}
D_X=\nabla_X+X \wedge \theta
\end{equation*}
for all $X$ in $TM$, where $\theta$ is in $\Lambda^1M$, in other words $\eta_X=X \wedge \theta$ for all $X$ in $TM$. A straightforward computation yields 
\begin{equation*}
[\eta_X, \eta_Y]-\eta_{Tor(X,Y)}=\vert \theta \vert^2 X \wedge Y
\end{equation*} 
for all $X,Y$ in $TM$, hence the comparison formula \eqref{Kcomp} becomes
\begin{equation*} \label{curvvect}
R^D(X,Y)=R(X,Y)+X \wedge D_Y \theta+D_X \theta \wedge Y+\vert \theta \vert^2 X \wedge Y
\end{equation*} 
whenever $X,Y$ in $TM$. Now an elementary computation shows that
\begin{pro} \label{Vvec}
Let $D$ belong to the class $W_3$, that is $\eta_X=X \wedge \theta$ for all $X$ in $TM$ and for some $\theta$ in $\Lambda^1M$. The following hold:
\begin{itemize}
\item[(i)] $(b_1R^D)_X=X \wedge d\theta $ for all $X$ in $TM$;
\item[(ii)] $ R^D(X,Y,Z,U)-R^D(Z,U,X,Y)=-\langle [X \wedge Y, d\theta], Z \wedge U\rangle $
for all $X,Y,Z,U$ in $TM$.
\end{itemize}
\end{pro}
Since this is needed to clarify when (iii) in Theorem \ref{res2} occurs we shall briefly examine 
the geometry of flat connections  with vectorial torsion. If $D$ is a metric connection with vectorial torsion determined by $\theta$ in $\Lambda^1M$ we form a linear connection on $M$ by $D^W=D-\theta \otimes 1_{TM}$. It is torsion 
free and satisfies $D^Wg=2\theta \otimes g$, hence it preserves the conformal class $c$ of $g$. In other words $D^W$ is a Weyl derivative and $(c,D^W)$ is a Weyl structure on $M$ (see \cite{CP} 
for more details). The curvature tensors of our two connections are related by $R^{D^{W}}=R^D-d\theta \otimes 1_{TM}$. If moreover $d\theta=0$ the curvature tensors coincide hence the holonomy algebras of 
$D^W$ and $D$ are equal by using the Ambrose-Singer holonomy theorem and then $Hol^0(D^W)=Hol^0(D)$ as well. Otherwise the Weyl structure and the vectorial torsion one are unrelated.
  
\begin{pro} \label{fvec}  A connection $D$ with vectorial torsion is flat if and only if its associated Weyl structure $(c,D^W)$ is flat. 
\end{pro}
\begin{proof}
If $D$ is flat $d\theta=0$ by using the algebraic Bianchi identity in (i) of Proposition \ref{Vvec}  whence the flatness of $D^W$. If the latter is flat $R^D=d\theta \otimes 1_{TM}$ hence $d\theta=0$ by using that $D$ is a metric connection.

\end{proof}
On compact manifolds, flat Weyl structures $(M^n,c,D^W)$ are well understood (see \cite{pgweyl}). Indeed, if $g_G$ is the Gauduchon metric in the conformal class $c$ then $(M,g_G)$ is locally isometric to $S^1 \times S^{n-1}$. It is not clear though 
if our flat connection with vectorial torsion can have irreducibly acting holonomy. As well known this does not occur in the case when $(M,g)$ is simply connected and complete.
\subsection{Prolongations of Lie algebras}
Let $(V^n,g)$ be a Euclidean vector space. Given a Lie algebra $\g$ we shall say that $(\g,V)$ is an orthogonal representation if there exists 
a Lie algebra morphism $\g \to \mathfrak{so}(V)$. In case this is injective $(\g,V)$ is called faithful and the Lie algebra $\g$ can be identified 
with a Lie subalgebra of $\mathfrak{so}(V)$.Unless otherwise specified we shall deal in what follows only with faithful orthogonal representations.

For any Lie sub-algebra $\g \subseteq \mathfrak{so}(V)$ we define its first skew-symmetric prolongation by 
\begin{equation}
\Lambda^3V \cap (\Lambda^1V \otimes \g)=\{T \in \Lambda^3V : x \lrcorner T \in \g \ \mbox{for \ all \ } x \in V\}
\end{equation}
When $\g=\{0\}$ or $\mathfrak{so}(V)$, the skew-symmetric prolongation space equals $\{0\}$ or $\Lambda^3V$ respectively, therefore we shall assume in 
what follows that $\g$ is proper. 
\begin{rema} If $V$ is a real vector space and $\g$ a Lie subalgebra of $End_{\mathbb{R}}V$ the first prolongation $\g^{(1)}$ of $\g$ is usually defined 
by $\g^{(1)}=\{\beta : V \to \g : \beta_xy=\beta_yx \ \mbox{for all} \ x,y \in V\}$. Clearly, if $\g$ is a subalgebra of $\mathfrak{so}(V)$ we must have $\g^{(1)}=\{0\}$.
\end{rema}
The skew-symmetric prolongation enters naturally in the sequence
$$0 \to \Lambda^3V \cap (\Lambda^1V \otimes \g) \to \Lambda^3V \stackrel{\varepsilon^{\perp}}{\to} \Lambda^1V \otimes \g^{\perp},
$$
where the map $\varepsilon^{\perp}$ is obtained by considering the splitting 
$$ id=\varepsilon+\varepsilon^{\perp}
$$
along $\Lambda^1V \otimes \Lambda^2V=(\Lambda^1V \otimes \g) \oplus (\Lambda^1V \otimes \g^{\perp})$. Explicitly, 
\begin{equation*}
\varepsilon^{\perp}(T)_X=(X \lrcorner T)_{\g^{\perp}}
\end{equation*}
for all $X$ in $V$ and whenever $T$ is in $\Lambda^3V$. By means of an orthogonality argument  it is easy to see that moreover
\begin{equation} \label{spin7g}
\Lambda^1V \otimes \g^{\perp}=Im \varepsilon^{\perp} \oplus \biggl [ (\Lambda^1V \otimes \g^{\perp}) \cap Ker(a) \biggr ].
\end{equation}
The projection map $\varepsilon^{\perp}$ plays an important r\^ole in the theory of metric connections in the class $W_1$. For, given a 
$G$-manifold $(M^n,g,G)$ such a $G$-connection exists iff the intrinsic torsion tensor $\eta$ belongs to $Im \varepsilon^{\perp}$ \cite{Fri1, Fri2}. 
\begin{rema}
Following \cite{Fri2} we note that those $G$-structures such that $\varepsilon^{\perp}$ is an isomorphism always admit a (unique) 
connection with totally skew-symmetric torsion.
Therefore it is of interest to characterise this maximal case. When $G$ is a Lie group having the involution property it has been showed in \cite{Fri2} that the 
only such representation is $(Spin(7), \mathbb{R}^8)$. 
\end{rema}
In fact, $\varepsilon^{\perp}$ is an isomorphism if and only if the $\g$-modules 
$$\Lambda^3V \cap (\Lambda^1V \otimes \g), (\Lambda^1V \cap \g^{\perp}) \cap Ker(a)$$ 
vanish identically. We shall use 
this observation later on to describe, in full generality, the case when $\varepsilon^{\perp}$ is an isomorphism. An important case when the skew-symmetric 
prolongation space is a priori understood is described below.
\begin{teo} \label{fixspin}\cite{agfr} Let $(\g,V)$ be an orthogonal representation and suppose that there exists a spinor $\psi \neq 0$ such that $\g \psi=0$. Then 
$\Lambda^3V \cap (\Lambda^1V \otimes \g)=\{0\}$.
\end{teo}
\subsection{Algebraic curvature tensors and Berger algebras}
In order to study the $\g$-module $\Lambda^3V \cap (\Lambda^1V \otimes \g)$ in the case when the defining 
representation $(\g,V)$ is irreducible, we need to review a number of facts related to formal 
curvature tensors and to the Berger algebra of $(\g,V)$. To begin with, the space of 
algebraic curvature tensors on $V$ is given by 
\begin{equation*}
\begin{split}
\K(\mathfrak{so}(V))=&\{R \in \Lambda^2V \otimes \Lambda^2V : b_1 R=0 \}\\
=& S^2(\Lambda^2V) \cap Ker(a)
\end{split}
\end{equation*}
where $b_1: \Lambda^2V \otimes \Lambda^2V \to \Lambda^1V \otimes \Lambda^3V$ is the Bianchi map defined in Section 2 and $a : S^2(\Lambda^2V)\to \Lambda^4V$ is the total alternation map.

For any Lie sub-algebra $\g \subseteq \mathfrak{so}(V)$ the space of 
$\g$-valued curvature tensors of $V$ is given by  
\begin{equation*}
\K(\g,V)=\K(\mathfrak{so}(n)) \cap (\Lambda^2V \otimes \g).
\end{equation*}
\begin{defi} \label{berg}
The Berger algebra $\underline{\g}$ of the the metric representation $(\g,V)$ is the smallest subspace $\mathfrak{p} \subseteq \g$ such that $\K(\mathfrak{p},V)=\K(\g,V)$. 
\end{defi}
Directly from this definition it follows that if $\mathfrak{h} \subseteq \g$ satisfies $\K(\mathfrak{h},V)=\K(\g,V)$, then $\underline{\g} \subseteq \mathfrak{h}$ and also 
that 
\begin{equation*}
\underline{\g}=span \{ R(X,Y) : X,Y \in V, R \in \mathcal{K}(\g,V) \}.
\end{equation*}
The case when $\underline{\g}=\g$ is of special interest and then $(\g,V)$ is called a Riemannian holonomy representation. We recall that the Berger list of 
irreducible, orthogonal representations is
\begin{equation} \label{bergerlist}
\begin{tabular}{ | r |c | c | } 
\hline
$\mathfrak{g}$  &  V &  \\ 
\hline 
$\mathfrak{so}(n) $ & $\mathbb{R}^n$ & \\
$\mathfrak{u}(m) $ & $\mathbb{R}^{2m}$ & \\
$\mathfrak{su}(m)$ & $\mathbb{R}^{2m}$ &\\
$\mathfrak{sp}(m) \oplus \mathfrak{sp}(1)$ & $\mathbb{R}^{4m}$ & \\
$\mathfrak{sp}(m)$ & $\mathbb{R}^{4m}$ & \\
$\mathfrak{spin}(7)$ & $\mathbb{R}^8$ & \\
$ \mathfrak{g}_2 $ & $\mathbb{R}^7$ & \\
\hline
\end{tabular}
\end{equation}
For irreducible orthogonal representations a precise description of the space of formal curvature tensors is given 
below (see also \cite{MS} for the non-metric case).
\begin{teo}\cite{cs} \label{clascurv}
Let $(\g,V)$ be an irreducible, orthogonal representation. If $(\g,V)$ is faithful then either:
\begin{itemize}
\item[(i)] $\K(\g,V)=\{0\}$;
\item[(ii)] $(\g,V)$ belongs to the Berger list;
\item[(iii)] $\K(\g,V) \cong \mathbb{R}$;
\item[] or 
\item[(iv)] $(\g,V)=(\mathfrak{sp}(m) \oplus \mathfrak{u}(1),[[EL]])$.
\end{itemize}
\end{teo}
When $\K(\g,V) \cong \mathbb{R}$ the Lie algebra $\g$ is the isotropy algebra of a symmetric space of compact type and the isomorphism is realised by the Ricci 
contraction map.
Also note that instances in (ii) and (iii)
correspond to the case when $(\g,V)$ is a Riemannian holonomy representation, that is $\underline{\g}=\g$. For a proof of these facts and more detail 
on the representations in the Berger list we refer the reader to \cite{br1, br2, cs}. We equally need to recall that the main result in \cite[Thm.4.6]{cs} also implies:
\begin{pro} \label{grpvalg}Let $(G,V)$ be an irreducible orthogonal representation of a Lie group $G$, with Lie algebra $\g$. If the representation $(\g,V)$ 
is not irreducible then $\mathcal{K}(\g,V)=\{0\}$. 
\end{pro}
We end this section with a few facts, the first of which is well known \cite{Hitchin}, concerning the Ricci tensors of elements in $\mathcal{K}(\g,V)$, when in presence 
of additional invariant objects.The second fact below has been proved in \cite{cs} for irreducible representations and here we give it a direct proof. 
\begin{pro} \label{ricci} Let $(\g,V^n)$ be a metric representation. The following hold:
\begin{itemize}
\item[(i)] if there exists a spinor $x \neq 0$ such that $\g x=0$ then $Ric(R)=0$ for all $R$ in $\mathcal{K}(\g,V)$;
\item[(ii)] the same conclusion holds if $(\g,V)$ is faithful and there exists a non-degenerate, $\g$-invariant element, in 
$\Lambda^3V \cap (\Lambda^1V \otimes \g^{\perp})$.
\end{itemize}
\end{pro}
\begin{proof}
(i) Let $R$ belong to $\mathcal{K}(\g,V)$, so that $R(X,Y)x=0$ for all $X,Y$ in $V$. Since the Clifford contraction 
$ \sum \limits_{i=1}^ne_iR(e_i,X)x=0$ for all $X$ in $V$, after using that $b_1R=0$ we arrive at $(RicX)x=0$ for all $X$ in $V$ and the claim 
follows.\\
(ii) Let $T$ in $\Lambda^3V \cap (\Lambda^1V \otimes \g^{\perp})$ be $\g$-invariant and non-degenerate. If $R$ belongs to $\mathcal{K}(\g,V)$ we have 
$[R(X,Y),T]=0$ for all $X,Y$ in $V$. In particular, 
$$ \sum \limits_{i=1}^n e_i \lrcorner [R(e_i,X),T]=0
$$
for all $X$ in $V$, where $\{e_i, 1 \le i \le n\}$ is some orthonormal basis. But
\begin{equation*}
\begin{split}
\sum \limits_{i=1}^n e_i \lrcorner [R(e_i,X),T]=&\sum \limits_{i,k} e_i \lrcorner (R(e_i,X)e_k \wedge e_k \lrcorner T)\\
=& -RicX \lrcorner T-\sum \limits_{i,k} R(e_i,X)e_k \wedge e_i \lrcorner e_k \lrcorner T\\
=& -RicX \lrcorner T+\frac{1}{2}\sum \limits_{i,k} R(e_k,e_i)X \wedge e_i \lrcorner e_k \lrcorner T
\end{split}
\end{equation*}
for all $X$ in $V$, after using the Bianchi identity for $R$. Now the last sum equals $\sum \limits_{p=1}^n R(e_p \lrcorner T)X \wedge e_p$ and since
$X \lrcorner T \in \g^{\perp}$ for all $X$ in $V$, it vanishes. We have shown that $RicX \lrcorner T=0$ for all $X$ in $V$, therefore 
the non-degeneracy of $T$ yields the vanishing of the Ricci contraction of $R$.
\end{proof}
\section{Structure results and linear holonomy algebras} \label{strres}
We shall compute in this section the first skew-symmetric prolongation of a proper Lie sub-algebra $\g $ of $\mathfrak{so}(V)$, i.e the intersection 
$\Lambda^3V \cap (\Lambda^1V \otimes \g)$, under the additional assumption 
that the representation $(\g,V)$ is irreducible. This will be done in several steps, by using the structure results on the space $\K(\g,V)$ in Theorem \ref{clascurv}.

If $\gamma$ belongs to $\Lambda^{\star}V$, let $L_{\gamma} : \Lambda^{\star}V \to \Lambda^{\star}V$ be the exterior multiplication with $\gamma$, and let 
$L_{\gamma}^{\star}$ be the adjoint of $L_{\gamma}$. Note that if $\alpha=g(F, \cdot, \cdot)$ is a $2$-form then one has 
$$ L_{\alpha}^{\star} \varphi=\frac{1}{2} \sum \limits_{i=1}^n Fe_i \lrcorner e_i \lrcorner \varphi
$$
whenever $\varphi$ belongs to $\Lambda^{\star}V$, where $\{e_i, 1 \le i \le n\}$ is some orthonormal basis in $V$. We start by collecting below a few simple observations on the tensors $R^T, \Omega^T$ as they have been 
defined in section 2.1 for a given $3$-form $T$.

\begin{lema} \label{curv1} Let $T$ belong to $\Lambda^3V \cap (\Lambda^1V \otimes \g)$. The following hold:
\begin{itemize}
\item[(i)] $R^T$ belongs to $\K(\g,V)$;
\item[(ii)] the $4$-form $\Omega^T$ belongs to $\Lambda^4V \cap (\Lambda^2V \otimes \g)$;
\item[(iii)]  $R^T=0$ if and only if $T=0$;
\item[(iv)] if moreover $(\g,V^n)$ is irreducible with $\mathcal{K}(\g,V) \cong \mathbb{R}$ we have 
\begin{equation*}
\langle X \lrcorner T, Y \lrcorner T \rangle=\frac{3}{n} \vert T \vert^2 \langle X,Y\rangle  
\end{equation*}
for all $X,Y$ in $V$.
\end{itemize}
\end{lema}
\begin{proof}
(i) and (ii) follow by using the definition of $R^T$, Lemma \ref{curvT1}, (i) together with \eqref{locomega} and the fact that $T_X$ belongs to $\g$ for all $X$ in $V$, by assumption.
(iii) if $R^T=0$ its Ricci contraction must vanish and the claim follows from Lemma \ref{curvT1}, (ii) by a positivity argument.
(iv) since $\mathcal{K}(\g,V)$ is $1$-dimensional any of its elements is $\g$-invariant. In particular $R^T$ is $\g$-invariant hence so is 
its Ricci contraction. Because our representation is irreducible, the claim follows now from Lemma \ref{curvT1}, (ii). 
\end{proof}
\begin{rema} \label{sign}
(i) Part (i) in the Lemma above still holds, with unchanged proof, if the metric $g$ has arbitrary signature.\\
(ii) If $g$ has indefinite signature, there are examples \cite{CorScha} of $3$-forms $T \neq 0$ such that $T_X T_Y=0$ for all $X,Y$ in $V$, and hence $R^T=0$.
\end{rema}
\begin{pro}\label{step1} Let $(\g,V)$ be irreducible and proper. The following hold:
\begin{itemize}
\item[(i)] if $\K(\g,V)=0$ we must have $\Lambda^3V \cap (\Lambda^1V \otimes \g)=\{0\}$;
\item[(ii)] if $\K(\g,V)$ is $1$-dimensional then either $\Lambda^3V \cap (\Lambda^1V \otimes \g)=\{ 0 \}$ or $(\g,V)$ is the adjoint representation of a compact, simple 
Lie algebra. In the latter case, $\Lambda^3V \cap (\Lambda^1V \otimes \g)$ is $1$-dimensional, 
generated by the Cartan $3$-form of $\g$.
\end{itemize}
\end{pro}
\begin{proof}
(i) follows directly from Lemma \ref{curv1}, (i) and (iii).\\
(ii) Let us assume that $T \neq 0$ belongs to $\Lambda^3V \cap (\Lambda^1V \otimes \g)$. 
If $\z$ in $V$ is fixed then $\z \lrcorner \Omega^T$ belongs to $\Lambda^3V \cap (\Lambda^1V \otimes \g)$ by Lemma \ref{curv1}, (ii) hence by (iv) in the same Lemma 
\begin{equation*}
\langle X \lrcorner \z \lrcorner \Omega^T,  Y \lrcorner \z \lrcorner \Omega^T \rangle=\frac{3}{n}\vert \z \lrcorner \Omega^T \vert^2 \langle X,Y\rangle
\end{equation*}
for all $X,Y$ in $V$. Taking $X=Y=\z$ leads easily to $\Omega^T=0$. By using \eqref{locomega} it follows that 
$R^T(X,Y)=3T_{T_XY}$ for all $X,Y$ in $V$. In other words $R^T \in \mathcal{K}(\mathfrak{h},V)$, where the subalgebra $\mathfrak{h} \subseteq \g$ is given by 
$$ \mathfrak{h}=\{T_X : X \in V\}.
$$
Since $0 \neq \mathcal{K}(\mathfrak{h},V) \subseteq \mathcal{K}(\g,V)$ and the latter is $1$-dimensional we conclude that  
$\mathcal{K}(\mathfrak{h},V)=\mathcal{K}(\g,V)$ hence $\underline{\g} \subseteq \mathfrak{h}$. Because 
$\underline{\g}=\g$ we must have $\mathfrak{h}=\g$, in other words $(\g,V)$ is the adjoint representation of a compact, simple Lie algebra (see also \cite{Rohr} for a related question).

In particular, any element $T^{\prime} \neq 0$ of $\Lambda ^3V \cap (\Lambda^1V \otimes \g)$ must be $\g$-invariant since it generates $\g$ as above and satisfies $\Omega^{T^{\prime}}=0$. But any such element 
can be written as $T^{\prime}_X=T_{QX}$ for all $X$ in $V$, where 
$Q$ is an endomorphism of $V$, hence $[\g,Q]=0$. The irreducibility of $(\g,V)$ makes that $Q=\lambda1_V+F$ where $\lambda$ is in $\mathbb{R}$ and $F$ is skew symmetric and $\g$-invariant. That is, 
$[T_X,F]=0$ for all $X$ in $V$ which leads, after total alternation, to $T_{FX}=0$ for all $X$ in $V$. The vanishing of $F$ is again granted by the irreducibility of $(\g,V)$ and 
the proof of the last claim is complete.
\end{proof}
\begin{rema} \label{ssign}
From the proof of (i) above we see that $\Lambda^3V \cap (\Lambda^1V \otimes \g)$ vanishes when $\mathcal{K}(\g,V)$ does, regardless of the irreducibility 
of the representation involved.
\end{rema}
We are now ready to prove our main result in this section.
\begin{teo} \label{class}
Let $\g \subseteq \mathfrak{so}(V^n)$ be proper such that $(\g,V)$ is irreducible. The following hold:
\begin{itemize}
\item[(i)] $\Lambda^3V \cap (\Lambda^{1}V \otimes \g)=\{0\}$ unless $(\g,V)$ is the adjoint representation of a compact, simple Lie algebra;
\item[(ii)] $\Lambda^pV \cap (\Lambda^{p-2}V \otimes \g)=\{ 0 \}$ for $4 \le p \le n$.
\end{itemize}
\end{teo}
\begin{proof}
(i) By Proposition \ref{step1} it remains to study the representations in the Berger list \eqref{bergerlist}. When 
$\g=\mathfrak{su}(n), \mathfrak{sp}(m),\mathfrak{sp}(m) \oplus \mathfrak{u}(1), \mathfrak{spin}(7), \mathfrak{g}_2$ the skew-symmetric prolongation vanishes,
because in these cases the algebra $\g$ fixes a spinor and Theorem \ref{fixspin} applies. When $\g=\mathfrak{u}(m)$ or $\mathfrak{sp}(m) \oplus \mathfrak{sp}(1), m \ge 2 $ it is 
an easy direct exercise to verify that the skew-symmetric prolongation vanishes.\\
(ii) the proof consists in the examination of the various instances for $\Lambda^3V \cap (\Lambda^{1}V \otimes \g)$, as indicated in (i).
\begin{itemize}
\item[(1)] $\Lambda^3V \cap (\Lambda^{1}V \otimes \g)=\{0\}$;\\
If $X_k, 1 \le k \le p-3$ in $V$ are fixed and $T$ belongs to $\Lambda^pV \cap (\Lambda^{p-2}V \otimes \g)$, then $X_1 \lrcorner \ldots X_{p-3} \lrcorner T $ belongs to $\Lambda^3V \cap (\Lambda^{1}V \otimes \g)$ and hence it 
vanishes, making that $T=0$.
\item[(2)] $\Lambda^3V \cap (\Lambda^{1}V \otimes \g) \neq 0$;\\
In analogy with the inductive argument above it is enough to see that $\Lambda^4V \cap (\Lambda^2V \otimes \g)=\{ 0 \}$. Since in this 
case $(\g, V)$ is the adjoint representation of a simple Lie algebra, let $\varphi$ in $\Lambda^3V \cap (\Lambda^{1}V \otimes \g)$ be its Cartan form. Then any 
element $\beta$ in $\Lambda^4V \cap (\Lambda^2V \otimes \g)$ satisfies 
$$ X \lrcorner \beta=\langle X, \z\rangle \varphi
$$
for all $X$ in $V$, where $\z$ belongs to $V$. If $\z \neq 0$, by using that $\beta$ is a $4$-form we get that $\z \lrcorner \varphi=0$, a contradiction since the Cartan form of a simple Lie algebra
is non-degenerate. Therefore $\z=0$, hence $\beta=0$, and we have showed that  $\Lambda^4V \cap (\Lambda^2V \otimes \g)=\{ 0 \}$ in this case as well. 
\end{itemize}
\end{proof}

Theorem \eqref{class} can now be directly applied to obtain a classification of the holonomy algebras generated by a $3$-form or equivalently of the holonomy algebras of 
connections with constant torsion in flat space. Indeed, given  a linear $3$-form $T \neq 0$ in $\Lambda^3V$ one forms the connection 
$$ \nabla_X^T=\nabla_X+\frac{1}{2}T_X
$$
for all $X$ in $V$. Here $V$ is equipped with its flat metric and $\nabla$ is the associated Levi-Civita connection. The connection $\nabla^T$ is metric and 
its torsion is totally skew-symmetric, given by the form $T$. The holonomy algebra $\mathfrak{h}^{\star}_T$ is computed as 
$\mathfrak{h}^{\star}_T=[\g^{\star}_T,\g^{\star}_T]$ where 
the Lie subalgebra $\g_T^{\star}$ of $\mathfrak{so}(V)$ is defined by 
$$ \g_T^{\star}=Lie \{T_X : X \in V \}.
$$
It is known \cite{agfr} that $\g_T^{\star}$ is semisimple hence $\mathfrak{h}^{\star}_T=\g^{\star}_T$. It has also been shown in \cite{agfr} that 
for any $T$ there is an orthogonal splitting of representations 
\begin{equation} \label{linsplit}
(\g_T^{\star},V)=V_0 \oplus \bigoplus \limits_{i=1}^d (\g^{\star}_{T_i},V_i)
\end{equation}
where $V_0$ is the trivial factor and the representations $(\g_{T_i}^{\star},V_i)$ where $T_i$ is in $\Lambda^3V_i$ are irreducible for $1 \le i \le d$.

In other words the study of the representation $(\g_T^{\star},V)$ reduces to the case when it is irreducible. 
$\\$
{\bf{Proof of Theorem \ref{hol3} }}:\\
Since by construction we have that 
$T$ belongs to $\Lambda^3V \cap (\Lambda^1V \otimes \g_T^{\star})$ it suffices to apply Theorem \ref{class}.
\subsection{Reducible representations}
To compute $\Lambda^3V \cap (\Lambda^1V \otimes \g)$ in full generality it is necessary to examine the case when the representation $(\g,V)$ is reducible. Indeed, we 
consider the orthogonal splitting 
$$V=V_0 \oplus \bigoplus \limits_{k=1}^d V_k$$
where $\g$ acts trivially on $V_0$ and each of the representations $(\g,V_k), 1 \le k\le d$ are irreducible. Moreover we define 
$$ \hat{V}_i=\bigoplus \limits_{\stackrel{k \neq i}{k=1}}^{d} V_k
$$
and let $\hat{\pi}_k : \g \to \mathfrak{so}(\hat{V}_k)$ the restriction of our 
original representation to $\hat{V}_k$. Finally, we consider the ideal  
$$ \hat{\g}_k:=(Ker \hat{\pi}_k)
$$
of $\g$ given by for all $1 \le k \le d$. Then each of the orthogonal representations $(\hat{\g}_k,V_k)$ is faithful.
\begin{pro} \label{split1}
Let $(\g,V)$ be faithful and orthogonal. We have that 
\begin{equation*}
\Lambda^3V \cap (\Lambda^1V \otimes \g)= \bigoplus \limits_{k=1}^d \Lambda^3V_k \cap (\Lambda^1V_k \otimes \hat{\g}_k).
\end{equation*}
\end{pro}
\begin{proof} Let us pick $T$ in $\Lambda^3V \cap (\Lambda^1V \otimes \g)$, that is $T$ in $\Lambda^3V$ such that $T_v$ belongs to $\g$, for all $v$ in $V$. 
First of all, since $\g V_0=0$ we get that $T_vV_0=0$ for all $v$ in $V$ hence $T_{v_0}=0$ for all $v_0$ in $V_0$. Now since $\g V_k \subseteq V_k$ it follows that 
$T_v V_k \subseteq V_k$ for all $v$ in $V$ and for all $1 \le k \le d$, in other words 
$ T \ \mbox{is in} \ \bigoplus \limits_{k=1}^d \Lambda^3V_k.$
Let us write $T=\sum \limits_{k=1}^d T_k$ for the corresponding decomposition of $T$. Then for 
any $v$ in $V_k, 1 \le k \le d$ we have $T_v\hat{V}_k=0$ and since $T_v$ belongs to 
$\g$ it follows that $T_k$ belongs to $\Lambda^3V_k \cap (\Lambda^1V_k \otimes \hat{\g}_k)$ for all $1 \le k \le d$ and the claim follows easily.
\end{proof}
The decomposition algorithm used above has been first presented and used to obtain a similar result for the space $\mathcal{K}(\g)$ in \cite{cs}. Since the 
representations $(\hat{\g}_i,V_i)$ may be still reducible, it is helpful to make the following straightforward:
\begin{lema} \label{redreps}Let $(\g,V)$ be a faithful orthogonal representation admitting a $\g$-invariant and orthogonal splitting $V=V_1 \oplus V_2$. If
$\mathfrak{l}$ is the ideal of $\g$ given by $\mathfrak{l}=(\hat{\g}_1 \oplus \hat{\g}_2)^{\perp}$ then:
\begin{itemize}
\item[(i)] the representations $(\hat{\g}_k \oplus \mathfrak{l},V_k), k=1,2$ are faithful;
\item[(ii)] if $\g$ acts irreducibly on $V_1$ then $(\hat{\g}_1 \oplus \mathfrak{l},V_1)$ is an irreducible representation. 
\end{itemize}
\end{lema}
Whereas this does not iterate well, it can be used  for the computation of a $\g$-module 
splitting as the skew-symmetric prolongation does.
\begin{teo} \label{gcase}
Let $(\g,V)$ be a faithful orthogonal representation. Then $(\g,V)$ splits as a direct sum of representations of the following type:
\begin{itemize}
\item[(i)] $(\mathfrak{so}(W),W)$ where $W$ is some Euclidean vector space;
\item[(ii)] adjoint representations of compact, simple Lie algebras;
\item[(iii)] representations with vanishing skew-symmetric prolongation.
\end{itemize}
\end{teo}
\begin{proof}
If $\Lambda^3V \cap (\Lambda^1V \otimes \g)=\{0\}$ the representation is of type (iii) and there is nothing to prove. We assume now that the skew-symmetric prolongation of $(\g,V)$ does not vanish and that 
our representation is reducible. Since $(\g,V)$ is faithful, it cannot be trivial, therefore 
Proposition \ref{split1} implies the existence of some $1 \le p \le d$ with $\Lambda^3V_p \cap (\Lambda^1V_p \otimes \hat{\g}_p) \neq \{0\}$.
We split orthogonally $\g=\hat{\g}_p \oplus \check{\g}_p \oplus \mathfrak{l}$ where $\check{\g}_p=\mathfrak{so}(\hat{V}_p) \cap \g $ and 
$\mathfrak{l}=(\hat{\g}_p \oplus \check{\g}_p)^{\perp}$ and note that $(\hat{\g}_p \oplus \mathfrak{l},V_p)$ is faithful and irreducible by Lemma \ref{redreps}.
Moreover $\Lambda^3V_p \cap (\Lambda^1V_p \otimes (\hat{\g}_p \oplus \mathfrak{l}) ) \neq \{0\}$ since it contains the non-vanishing  
$\Lambda^3V_p \cap (\Lambda^1V_p \otimes \hat{\g}_p)$.
By applying Theorem \ref{class}, (i) it follows that $\hat{\g}_p \oplus \mathfrak{l}$ is a simple Lie algebra. Because $\hat{\g}_p \neq 0$ either the ideal $\mathfrak{l}$ vanishes and 
hence $(\hat{\g}_p,V_p)$ splits out as a direct factor 
of type (i) or (ii), or $V_p$ is $4$-dimensional and $\hat{\g}_p=\Lambda^2_{+}V_p \cong \mathfrak{so}(3)$. But this last case is not eligible 
since then $(\hat{\g}_p,V_p)$ would be irreducible hence  with vanishing skew prolongation, again by Theorem \ref{class}, (i) and dimension comparison, a contradiction. The proof is now completed by an induction 
argument.
\end{proof}
\section{Uniqueness of connections with totally skew-symmetric torsion}
Let $(M^n,g,G)$ be a connected Riemannian $G$-manifold and associated representation $(G,V)$. If $\g$ denotes the Lie algebra of the group $G$ this gives rise to a subbundle of $\Lambda^2M$ with 
fibers isomorphic to $\g$, to be denoted by the same symbol. 

Let now $D$ be a metric connection in the class $W_1$, that is $D=\nabla+\frac{1}{2}T$ where $T$ belongs to 
$\Lambda^3M$. The holonomy group $Hol(D)$ of the connection $D$ is a Lie subgroup of $O(n)$. We assume now that $D$ preserves the given $G$-structure $Hol(D) \subseteq G$ hence the subbundle $\g$ of $\Lambda^2M$ must 
be parallel w.r.t. the connection $D$. 
\begin{rema} 
Conversely, if $D$ is a connection in the class $W_1$ let $G=Hol(D)$ be its holonomy group at some point of $M$, with Lie algebra $\g$. Since $\g$ is invariant under the action of $Hol(D)$ by using parallel transport we obtain a 
subbundle $\g$ of $\Lambda^2M$ which is preserved by the connection $D$. However the bundles obtained by this construction might be different when starting from two connections having the same holonomy.
\end{rema}
If there exists another metric connection 
$$D^{\prime}=\nabla+\frac{1}{2}T^{\prime}$$ 
in the class $W_1$ preserving the $G$-structure the subbundle $\g \subseteq \Lambda^2M$ is preserved by both $D$ and $D^{\prime}$. It follows that 
\begin{equation*}
[X \lrcorner (T-T^{\prime}),\g] \subseteq \g
\end{equation*}
for all $X$ in $TM$. In other words the difference $T-T^{\prime}$ belongs, at any point of $M$, to $\Lambda^3M \cap (\Lambda^1M \otimes \tilde{\g})$ where we have 
defined the extension $\tilde{\g}$ of the holonomy algebra $\g$ by 
\begin{equation*}
\tilde{\g}=\{ F \in \Lambda^2M : [F, \g] \subseteq \g\}.
\end{equation*}
Clearly, both $\Lambda^3M \cap (\Lambda^1M \otimes \tilde{\g})$ and $\tilde{\g}$ define subbundles of $\Lambda^3M$ and $\Lambda^2M$ respectively, since they parallel w.r.t. the connection $D$. \\
$\\$
{\bf{Proof of Theorem \ref{uniq1}:}}\\
We start working on the representation $(G,V)$ at a given point of the manifold $M$. Let us consider the Lie subgroup of $O(n)$ given by $\tilde{G}=\{g \in O(n) : g^{\star}\g \subseteq \g \}$. The Lie algebra of $\tilde{G}$ clearly equals 
$\tilde{\g}$ and moreover $(\tilde{G},V)$ is irreducible since $G \subseteq \tilde{G}$. As indicated above, we need to examine the skew-symmetric prolongation space 
$\Lambda^3V \cap (\Lambda^1V \otimes \tilde{\g})$. If $(\tilde{\g},V)$ does not act irreducibly then $\mathcal{K}(\tilde{\g})=\{0\}$ by Proposition \ref{grpvalg} and then 
the skew-symmetric prolongation space vanishes as well (see also Remark \ref{ssign}). When $(\tilde{\g},V)$ is irreducible, using 
Theorem \ref{class} yields that $\Lambda^3V \cap (\Lambda^1V \otimes \tilde{\g})=\{ 0 \}$ unless 
$(\tilde{\g},V)$ is the adjoint representation of a compact, simple, Lie algebra. If this is the case, $\tilde{\g}=\g$ hence the $D$-parallel bundle $\Lambda^3M \cap (\Lambda^1M \otimes \tilde{\g})$ has real 
rank one. Therefore, if $t$ of unit length is a generator of the latter subbundle it must satisfy $Dt=0$ and $\Omega^t=0$. That the last equation is equivalent to $t^2=1$ in the Clifford algebra bundle $Cl(M)$ has been observed in 
\cite{witt}.
\section{Killing frames on Riemannian manifolds}
Let $(M^n,g)$ be a Riemannian manifold and let $D$ be a metric connection in the class $W_1$ with torsion form $T$ in $\Lambda^3M$. We wish to investigate here which are the situations when 
the connection $D$ is flat. A key r\^ole in our approach will be played by the cubic deformation $D^{\frac{1}{3}}$ the use of which shall be also useful in arriving at a Lie algebra prolongation problem.
\begin{lema} \label{flat3form}
Let us assume that the connection $D$ is flat. The following hold:
\begin{itemize}
\item[(i)] $D^{\frac{1}{3}}T=0$;
\item[(ii)] $R^{\frac{1}{3}}(X,Y)=\frac{1}{9}([T_X,T_Y]+T_{T_XY})$ for all $X,Y$ in $TM$.
\end{itemize}
\end{lema}
\begin{proof}
(i) Since $D$ is flat, that is $R^D=0$ we obtain from \eqref{curv3form} that 
$$0=(R-\frac{1}{12}R^T)-\frac{1}{2} \Theta-(\frac{1}{4}d_DT+\frac{1}{3} \Omega^T).
$$
By identifying components along the splitting 
$$\Lambda^2M \otimes \Lambda^2M=\mathcal{K}(M) \oplus \Lambda^2(\Lambda^2M) \oplus \Lambda^4M$$ 
we obtain 
\begin{equation*}
\begin{split}
&R=\frac{1}{12}R^T,\\
&\Theta=0,\\
&d_DT=-\frac{4}{3}\Omega^T.
\end{split}
\end{equation*}
The vanishing of $\Theta$ combined with (i) in Lemma \ref{bianchis} yields 
$$ D_XT=\frac{1}{4}X \lrcorner d_DT=-\frac{1}{3}X \lrcorner \Omega^T=\frac{1}{3}[T_X,T]
$$
for all $X$ in $TM$. The parallelism of $T$ w.r.t $D^{\frac{1}{3}}$ follows now from the definition of the family of connections 
$D^t, t$ in $\mathbb{R}$.\\
(ii) Again from the vanishing of $R^D$ and $\Theta$ combined with the fact that  
$d_DT=-\frac{4}{3}\Omega^T$ it follows that \eqref{param} applied for the value 
$t=\frac{1}{3}$ actually reads 
\begin{equation*} \label{curvcube}
\begin{split}
R^{\frac{1}{3}}(X,Y)=&\frac{1-t}{12}((t+1)R^T-2t\Omega^T)(X,Y)\\
=&\frac{1}{9}([T_X,T_Y]+T_{T_XY})
\end{split}
\end{equation*}
for all $X,Y$ in $TM$,after using also the definition of $R^T$ and \eqref{locomega}. 
\end{proof}
We recall that 
the action $\mathfrak{so}(TM) \times \mathcal{K}(M) \to \mathcal{K}(M)$ is given by 
$$ [F,Q](X,Y)=Q(FX,Y)+Q(X,FY)+[Q(X,Y),F] $$
for all $X,Y$ in $TM$ and whenever $(F,Q)$ belongs to $\mathfrak{so}(TM) \times \mathcal{K}(M)$, where the usual identifications apply.\\
$\\$
{\bf{Proof of Theorem \ref{class3flat}}}:\\
From Lemma \ref{flat3form}, (i) it follows that $D^{\frac{1}{3}}R^T=0$. At the same time it has been shown in \cite{datri} that $(M^n,g)$ is a locally symmetric space that is $\nabla R=0$. Since $R$ is proportional to 
$R^T$ combining these two facts yields 
\begin{equation} \label{prol3form}
[T_X,R^T]=0
\end{equation}
for all $X$ in $TM$. In other words, at each point of $M$ the tensor $R^T$ is $\g^{\star}_T$-invariant. At a fixed point $m$ of $M$ we consider the splitting of 
\begin{equation} \label{split3f}
T_mM=V_0 \oplus \bigoplus \limits_{i=1}^d V_i
\end{equation}
into irreducible components under the action of $\g^{\star}_{T_m}$, where $V_0$ is a trivial factor.  By (ii) in Lemma \ref{flat3form}  the holonomy algebra $\mathfrak{hol}(D^{\frac{1}{3}}) \subseteq \g^{\star}_{T_m}$ hence 
the splitting is $\mathfrak{hol}(D^{\frac{1}{3}})$-invariant. Since $M$ is simply connected it is also $Hol(D^{\frac{1}{3}})$-invariant and hence extends by 
parallel transport to a $D^{\frac{1}{3}}$-parallel splitting of $TM$ to be denoted similarly. From \eqref{linsplit} it follows that $T$ is split along \eqref{split3f}, that is 
$$T=\sum \limits_{i=1}^d T_i
$$
where $T_i$ belongs to $\Lambda^3V_i, 1 \le i \le d$. It follows that each $V_i, 0 \le i \le d$ is parallel w.r.t. the Levi-Civita connection of $g$ 
hence by using the deRham splitting theorem for the complete metric $g$ we have that 
$(M,g)=(M_0,g_0) \times (M_1,g_1) \times \ldots (M_d,g_d)$, where $M_0$ is flat and $\g^{\star}_{T_i}$ acts irreducibly on each tangent space to $M_i, 1 \le i \le d$. 
Moreover the non-flat factors $M_i, 1\le i \le d$ fall, by Theorem \ref{hol3} into two classes: 
\begin{itemize}
\item[(i)] when  $\g^{\star}_{T_i}=\mathfrak{so}(TM_i)$. In this case it follows from \eqref{prol3form} that $R^{T_i}$ must have constant non-zero sectional curvature, hence the same is true for the metric $g_i$. 
Now it has been shown in \cite{datri} that the only possibilities here are up to homothety $S^3,S^7$ with the round metric.
\item[(ii)] when $(\g^{\star}_{T_i},TM_i)$ is at each point of $M_i$ the adjoint representation of a simple, compact Lie algebra. Then $\Omega^{T_i}=0$ vanishes hence $T_i$ is parallel w.r.t. the Levi-Civita connection of 
$g_i$ by using (i) in Lemma \ref{flat3form}. The fact that $(M_i,g_i)$ is isometric to a Lie group with a bi-invariant metric is now straightforward, see e.g. \cite{datri} for details.
\end{itemize}
\section{A class of Pl\"ucker type embeddings and $k$-Lie algebras} 
Using the computation of the skew-symmetric prolongation of a Lie algebra of compact type we have obtained previously, we shall present here 
a class of Pl\"ucker type relations, which have been first considered and treated in dimensions up to $8$ in \cite{ofarill}. Our present context will be that 
of a Euclidean vector space $(V^n,g)$. We are interested here in the class of forms $T$ in $\Lambda^pV, p \ge 3$ satisfying 
\begin{equation} \label{plc}
[L_{\z}^{\star}T,T]=0
\end{equation}
for all $\z$ in $\Lambda^{p-2}V$, where we recall that for any $\z$ in $\Lambda^{p-2}V$ the map $L^{\star}_{\z}$ is the adjoint of the exterior multiplication by $\z$.
When $p=3$ every form $T$ subject to \eqref{plc} induces a Lie algebra structure on $V$ which is given by $[X,Y]=T_XY$ 
for all $X,Y$ in $V$(see also \cite{Rohr} for a related question). We recall that a $p$-form is called {\it{decomposable}} if it can be written as the sum of mutually orthogonal simple 
$p$-forms. 

Note that the classical Pl\"ucker relations (see \cite{GrHa, Eastwodd}) for a $p$-form $T$ state that $T$ is decomposable if and only if it satisfies 
\begin{equation} \label{classic}
L^{\star}_{\z}T \wedge T=0
\end{equation}
whenever $\z$ belongs to $\Lambda^{p-1}V$. 

The Euclidean version of a conjecture in \cite{ofarill} states that any $p$-form, $p \ge 4$ satisfying \eqref{plc} on a Euclidean vector space must be decomposable. To prove 
this is indeed the case let us first observe the impact of having \eqref{plc} satisfied for some form $T$ on its isotropy algebra. Recall that the latter is 
defined by
$$\g_T=\{\alpha \in \Lambda^2V:[\alpha,T]=0\}$$ 
and also that the representation $(\g_T,V)$ is orthogonal and faithful. At the same time let us consider the Lie algebra
\begin{equation*}
\mathfrak{r}_T:=Lie \{L^{\star}_{\z}T: \z \in \Lambda^{p-2}V\}.
\end{equation*}
A form $T$ in $\Lambda^pV$ satisfies the Pl\"ucker relations \eqref{plc} if and only if one has $\mathfrak{r}_T \subseteq \g_T$. Both $\g_T$ and 
$\mathfrak{r}_T$ have good splitting properties and in order to progress in that direction we make:
\begin{defi} \label{ired}
Let $T$ belong to $\Lambda^pV$. We call $T$ irreducible iff the orthogonal representation $(\mathfrak{r}_T, V)$ is irreducible.
\end{defi}   
\begin{pro} \label{iredsplit}Let $T$ in $\Lambda^pV$ satisfy \eqref{plc}. Then $V$ admits an orthogonal, direct sum decomposition 
$$ V=V_0 \oplus \bigoplus \limits_{k=1}^{d}V_k
$$
such that $T$ vanishes on $V_0$ and $T=\sum \limits_{k=1}^{d}T_k$ where:
\begin{itemize}
\item[(i)] $T_k$ belong to $\Lambda^pV_k$ and satisfy \eqref{plc} for  $1 \le k \le d$;
\item[(ii)] the forms $T_k, 1 \le k \le d$ are irreducible.
\end{itemize}
\end{pro}
\begin{proof}
(i) since $(\mathfrak{r}_T,V)$ is an orthogonal representation it splits as 
$$ V=V_0 \oplus \bigoplus \limits_{k=1}^{d}V_k
$$
on orthogonal direct sum, where $\mathfrak{r}_T$ acts trivially on $V$ and $(\mathfrak{r}_T,V_k), 1 \le k \le d$ are irreducible. It follows that 
$\langle v_i \hook v_j \hook T, \z\rangle=0$ 
for all $v_i$ in $V_i, v_j$ in $V_j$ and for all $\z$ in $\Lambda^{p-2}V$. It is now easy to conclude that $T$ splits as indicated and (i) holds.\\
To prove (ii) we notice that since \eqref{plc} holds for $T_k$ we have that $\mathfrak{r}_{T_k} \subseteq \g_{T_k}$ for all $1 \le k \le d$, hence the latter 
acts irreducibly on $V_k$ as the former does. 
\end{proof}
$\\$
{\bf{Proof of Theorem \ref{plckclass}:}}\\
The defining equation \eqref{plc} actually says that $T$ belongs to 
$$(\Lambda^{p-2}V \otimes \g_T) \cap \Lambda^pV.$$ 
Let us assume now that $T$ is irreducible. Since $T \neq 0$, by using Theorem \ref{class}, (ii) it follows that $\g_T=\mathfrak{so}(V)$. In this case $T$ must be invariant under $\mathfrak{so}(V)$, hence 
proportional to a volume form. In particular this forces the equality $n=p$.\\
The general case follows now by using the splitting result in Proposition \ref{iredsplit}, (ii).
\subsection{Metric $n$-Lie algebras} 
The Pl\"ucker relations \eqref{plc} are related, as it has been observed in \cite{ofarill} to the class of the so-called $n$-Lie algebras. More precisely, 
\begin{defi}
Let $V$ be a vector space over $\mathbb{R}$ (not necessarily of finite dimension). Then a $n$-Lie algebra structure on $V$ consists in a map 
$$ [\cdot, \dots, \cdot] : \Lambda^nV \to V $$ 
such that the generalised Jacobi identity 
\begin{equation*}
[X_1, \ldots ,X_{n-1}, [Y_1, \dots, Y_n]]=\sum \limits_{i=1}^n [Y_1, \dots, [X_1, \ldots, X_{n-1}, Y_i], \ldots, Y_n]
\end{equation*}
holds, whenever $X_k,Y_k, 1 \le k \leq n-1$ are in $V$.\\
An $n$-Lie algebra is called metric if there is a non-degenerate, symmetric bilinear form $\langle \cdot , \cdot \rangle$ on $V$ such that 
$$ \langle [X_1, \ldots, X_n], X_{n+1}\rangle=-\langle [X_1, \ldots, X_{n+1}], X_{n} \rangle $$  
\end{defi}
Finally, let us call a metric $n$-Lie algebra Euclidean iff its associated bilinear form is positive definite. To see that one can classify 
finite dimensional, Euclidean $n$-Lie algebra we need to recall first 
\begin{pro} \cite{ofarill}
Let $(V, <\cdot, \cdot>, [\cdot, \ldots, \cdot])$ be a metric $n$-Lie algebra. Then the $n+1$-form $T$ in $\Lambda^{n+1}V$ defined by 
$$ T(X_1, \ldots, X_n, X_{n+1})=\langle [X_1, \ldots, X_n], X_{n+1} \rangle 
$$
whenever $X_k, 1 \le k \le n+1$ belong to $V$, satisfies the Pl\"ucker type relations \eqref{plc}. 
\end{pro}
This is very similar to the construction of the Cartan 
$3$-form of a metric Lie algebra. By using Theorem \ref{plckclass} we immediately obtain the following.
\begin{teo} \label{nLieclass}
Let $V$ be a finite dimensional Euclidean $n$-algebra. If the form $T$ is irreducible then either $V$ is a Lie algebra or $n=dim_{\mathbb{R}}V-1$, in which case 
$T$ is proportional to a volume form.
\end{teo}
\section{The holonomy of connections with vectorial torsion}
This section is devoted to the study of the holonomy representation of metric connections with torsion of vectorial type.
Particular attention will be paid to the properties of invariant $4$-forms, for those shall be used to obtain geometric information 
in our situation. 
\subsection{Invariant forms and Casimir operators}
Let us consider again a Euclidean vector space $(V^n,g), n \ge 5$ equipped with a faithful action of some Lie algebra $\g$, therefore considered as a 
subalgebra of $\mathfrak{so}(V)$. The existence problem of $\g$-invariant forms on $V$ has been considered in \cite{kos} by means of the following construction.
\begin{defi} \label{def4inv}
The characteristic form of the representation $(\g,V)$ is given by 
$$ T^{\g}=a(1_{\g})
$$
where 
$a : S^2\g \to \Lambda^4V$ is the alternation map.
\end{defi}
Clearly $T^{\g}$ in $\Lambda^4V$ is $\g$-invariant and if the representation $(\g,V)$ is induced by that of a Lie group $G$ then $T^{\g}$ is 
invariant under $G$ as well.To consider the case when the characteristic form vanishes we also need:
\begin{defi} \label{minrep}
Let $(\g,V)$ be a faithful metric representation. It is called minimal if one has $T^{\g}=0$.
\end{defi}
On a given representation minimality turns out to be a strong requirement as the following shows.
\begin{pro} \label{minrepcara}
Let $(\g,V)$ be faithful and metric. Then:
\begin{itemize}
\item[(i)] $(\g,V)$ is minimal iff $1_{\g}$ belongs to $\mathcal{K}(\g,V)$;
\item[(ii)] if there exists a spinor $x \neq 0$ such that $\g x=0$ then $T^{\g} \neq 0$.
\end{itemize}
\end{pro}
\begin{proof}
(i) follows immediately from the definition of $T^{\g}$, since $1_{\g}$ belongs to $S^2(\g)$.\\
(ii) supposing that $T^{\g}=0$, it follows by (i) that $1_{\g} \in \mathcal{K}(\g,V)$ hence $Ric(1_{\g})=0$ by Proposition \ref{ricci}, (i), since 
$\g$ preserves a spinor. But 
\begin{equation*}
\begin{split}
Ric(1_{\g})(X,Y)=&\sum \limits_{i=1}^n \langle (e_i \wedge X)_{\g}, e_i \wedge Y\rangle \\
=& \sum \limits_{x_k \in \g}\langle x_kX, x_kY \rangle 
\end{split}
\end{equation*}
for all $X,Y$ in $V$, where $\{x_k\}$ is an orthonormal basis in $\g \subseteq \mathfrak{so}(V) \cong \Lambda^2V $, w.r.t. the form inner product. It follows that $Ric(1_{\g})$ vanishes if and only if $(\g,V)$ is the trivial 
representation, a contradiction.
\end{proof}
All entries in the Berger list \eqref{bergerlist} but the first and  $\mathfrak{u}(m), \mathfrak{sp}(m) \oplus \mathfrak{sp}(1)$ fix a spinor and therefore 
cannot be minimal. Moreover a direct verification using this time (i) in the Proposition above shows that 
$\mathfrak{u}(m)$ and $\mathfrak{sp}(m) \oplus \mathfrak{sp}(1)$ are not minimal neither, making that only $\mathfrak{so}(n), n \neq 4$ is. On the contrary, 
the adjoint representation of a semisimple Lie algebra is minimal as well as any irreducible representation of a Lie algebra on a $5$-dimensional Euclidean space.
In fact, the list of minimal irreducible representations has been given in \cite{kos} in terms of their associated symmetric spaces.   

Various properties of a representation can be described in terms of the Casimir operators 
$$ C^k : \Lambda^kV \to \Lambda^kV
$$
of the representation $(\g,V)$, defined 
by $C^k(\varphi)=-\sum \limits_{x_k \in \g}[x_k, [x_k, \varphi]]$ for some orthonormal (w.r.t. the form inner product) basis $\{x_k \}$ in $\g  \subseteq \mathfrak{so}(V) \cong \Lambda^2V $ and for all 
$\varphi$ in $\Lambda^kV, 0 \le k \le n$. We have that $C^k \ge 0$ for all $k \ge 0$ and when 
$k=1$ the operator $C^1$ gives actually an element in $S^2V$, which is invariant under $\g$, that is $[\g,C^1]=0$. Therefore, when $(\g,V)$ is irreducible we must have $C^1=\mu1_V$ for some $\mu$ in 
$\mathbb{R}$ which is computed from the trace of $C^1$ by 
\begin{equation} \label{const}
\mu=\frac{2dim_{\mathbb{R}}\g}{n}.
\end{equation}

A technical observation we shall need is that a straightforward verification yields 
\begin{equation} \label{formal}
L^{\star}_{\alpha_1}(\alpha_2 \wedge \alpha_3)=g(F_3F_1F_2+F_2F_1F_3, \cdot)+\langle \alpha_1,\alpha_2\rangle \alpha_3+\langle \alpha_3,\alpha_1 \rangle \alpha_2 
\end{equation}
for all $\alpha_i=g(F_i \cdot, \cdot), 1 \le i \le 3$ in $\Lambda^2V$. It implies that for any orthogonal representation $(\g,V)$ we have
\begin{equation} \label{carca}
L^{\star}_F T^{\g}=2\pi_{\g}F+C^2F-\{C^1,F\}
\end{equation}
for all $F$ in $\Lambda^2V$, where $\pi_{\g} : \Lambda^2V \to \g$ is the orthogonal projection.
We now need two preparatory Lemmas which are directly related to 
the geometry of $G$-structures of vectorial type.
\begin{lema} \label{extr}
Let $(\g,V^n), n > 4$ be a faithful and irreducible orthogonal representation such that $\g \neq \mathfrak{so}(V)$. If $F$ in $\Lambda^2V$ satisfies 
\begin{equation} \label{vectoreq}
\begin{split}
&[F, \g^{\perp}] \subseteq \g \\
& L^{\star}_{F}T^{\g}=0
\end{split}
\end{equation}
then $F$ belongs to $\g^{\perp}$. 
\end{lema}
\begin{proof}
We split $F=A+B$ along $\Lambda^2V=\g \oplus \g^{\perp}$ and notice that 
$$[C^2F, \g^{\perp}]=[C^2A+C^2B, \g^{\perp}] \subseteq \g 
$$
because of $[\g,\g^{\perp}] \subseteq \g^{\perp}$ and after use of the Jacobi identity. Since our representation is irreducible, we have $C^1=\mu 1_V$ where $\mu$ is given by \eqref{const}.
Because $L^{\star}_FT^{\g}=0$, from \eqref{carca} we get after identifying components along $\Lambda^2V=\g \oplus \g^{\perp}$ that
\begin{equation*} \label{casidouble}
\begin{split}
&C^2A=(2\mu-2)A, \ C^2B=2\mu B
\end{split}
\end{equation*}
in particular $[(2\mu-2)A+2\mu B,\g^{\perp}] \subseteq \g$. Combined with our original assumption it yields $[A,\g^{\perp}] \subseteq \g$. Hence $[A,\g^{\perp}]=0$ by using again that 
$[\g,\g^{\perp}] \subseteq \g^{\perp}$. It follows that $\mathcal{C}A-C^2A=0$ where $\mathcal{C}$ is the Casimir operator of $\mathfrak{so}(V)$, which, as well known, is given by $\mathcal{C}=2(n-2)$.
It follows that $(2\mu-2-2(n-2))A=0$. If $A \neq 0$ by using \eqref{const} we get $dim_{\mathbb{R}} \g=\frac{n(n-1)}{2}$ whence $\g=\mathfrak{so}(V)$, a contradiction. Thus $A=0$ and the claim is proved.
\end{proof}
It is now necessary to supplement Lemma \ref{extr} by the following elementary
\begin{lema}  \label{symspace}
Let $(\g,V^n), n >4$ be a an orthogonal and faithful representation. 
\begin{itemize}
\item[(i)]If the space $\{F \in \g^{\perp} : [F, \g^{\perp}] \subseteq \g\}$ does not vanish then $\g$ is a symmetric Lie subalgebra of $\mathfrak{so}(V)$ such that $[\g^{\perp}, \g^{\perp}]=\g$;
\item[(ii)] if moreover $(\g,V)$ is irreducible and $F$ in $\Lambda^2V$ satisfies \eqref{vectoreq} then $F=0$.
\end{itemize}
\end{lema}
\begin{proof}
(i) Since the argument is standard it will only be outlined. Let $\mathfrak{m}$ be the space above which we assume to be non-zero. By definition $[\mathfrak{m}, \g^{\perp}] \subseteq \g$ hence 
\begin{equation*}
\begin{split}
[\mathfrak{m}, \mathfrak{m}] \subseteq \g,  [\mathfrak{m}, \mathfrak{m}^{\perp}] \subseteq \g.
\end{split}
\end{equation*}
Since $[\g, \mathfrak{m}] \subseteq \mathfrak{m}$ it follows that $[\mathfrak{m}, \mathfrak{m}^{\perp}]=0$. Define now the ideal $\mathfrak{i}_1$ in $\g$ by $\mathfrak{i}_1=[\mathfrak{m}, \mathfrak{m}]$ and let 
$\mathfrak{i}_2$ be the orthogonal complement of $\mathfrak{i}_1$ in $\g$. An orthogonality argument shows that $[\mathfrak{i}_2, \mathfrak{m}]=0$. Starting from $[\mathfrak{m}, \mathfrak{m}^{\perp}]=0$ the Jacobi identity yields 
$[[\mathfrak{m}, \mathfrak{m}], \mathfrak{m}^{\perp}]=0$ hence $[\mathfrak{i}_1, \mathfrak{m}^{\perp}]=0$. Therefore $[\mathfrak{m}^{\perp}, \mathfrak{m}^{\perp}]$ is orthogonal to $\mathfrak{i}_1 \oplus \mathfrak{m}$ and hence 
contained in $\mathfrak{i}_2 \oplus \mathfrak{m}^{\perp}$. It is now easy to see that $\mathfrak{i}_1 \oplus \mathfrak{m}$ as well as $\mathfrak{i}_2 \oplus \mathfrak{m}^{\perp}$ are subalgebras of $\mathfrak{so}(V)$, in fact ideals 
since the commutation rules displayed above lead to $[\mathfrak{i}_1 \oplus \mathfrak{m},\mathfrak{i}_2 \oplus \mathfrak{m}^{\perp}]=0$. We conclude by using that $\mathfrak{so}(V)$ is a simple Lie algebra as we have assumed that $n>4$.

(ii) Supposing that $F \neq 0$ it follows from Lemma \ref{extr} that $F$ belongs to $\g^{\perp}$ hence $\g$ is a symmetric Lie subalgebra of $\mathfrak{so}(V)$ by (i). Since $(\g,V)$ is irreducible, consulting the tables in \cite{Helgy} leads 
to $\g=\mathfrak{u}(\frac{n}{2})$ where $2n=dim_{\mathbb{R}}V$. In this case $T^{\g}$ is proportional to $\omega \wedge \omega$ where $\omega$ is the standard K\"ahler form on $\mathbb{R}^{2n}$. However, from \eqref{formal}
it follows that $L^{\star}_{\alpha}(\omega \wedge \omega)=2\alpha$ for all $\alpha$ in $\mathfrak{u}^{\perp}(\frac{n}{2})$ hence the second equation in \eqref{vectoreq} leads to $F=0$, a contradiction.
\end{proof}
\subsection{The classification}
The last ingredient we need in order to state and prove our main result in this section is the exactness criterion for closed $1$-forms in the following Lemma. 
\begin{lema} \label{lvglob} \cite{cs2}
Let $M$ be a connected manifold and let $\theta$ in $\Lambda^1M$ be a closed $1$-form, that is $d\theta=0$. If $\lambda$ is a smooth, non-identically zero, function on 
$M$ such that 
$$ d\lambda+\lambda \theta=0
$$ 
then $\lambda$ is nowhere vanishing and $\theta=-dln\vert \lambda \vert$.
\end{lema}
We can now make the following.
\begin{teo} \label{vec}
Let $(M^n,g), n \neq 4$ be a connected
and oriented Riemannian manifold and let $D$ be a metric connection with vectorial torsion 
$\theta$ in $\Lambda^1M \subseteq \Lambda^1M \otimes \Lambda^2M$. If the holonomy 
representation $(G,V)$ of $D$ is irreducible and proper then:
\begin{itemize}
\item[(i)] $d\theta=0$; 
\item[(ii)] if $(G,V)$ is not on the Berger list then either:
\begin{itemize}
\item[(a)] the connection $D$ is flat, 
\item[] or
\item[(b)] the metric $g$ is conformal to a non-flat,locally symmetric Riemannian metric;
\end{itemize}
\item[(iii)] if $(G,V)$ is on the Berger list then $(M,g)$ has a l.c.p. structure. 
\end{itemize} 
\end{teo}
\begin{proof}
Let us denote by $\g$ the Lie algebra of the Lie group $G$. At a given point of the manifold we consider the $G$-invariant form $T^{\g}$ which can be then extended over $M$, by using parallel transport, 
to a parallel form w.r.t to $D$. Since the curvature tensor $R^D$ of the connection $D$ lies 
in $\Lambda^2M \otimes \g$, by evaluation of (ii) in Proposition \ref{vec} on $\g^{\perp} \otimes \g^{\perp}$ we find that 
$$[d\theta, \g^{\perp}] \subseteq \g.$$ 
Since $DT^{\g}=0$ we have \cite{agfrivec} 
\begin{equation} \label{strveceqs}
dT^{\g}= 4\theta \wedge T^{\g}, \ d^{\star}T^{\g}=-(n-4)\theta \lrcorner T^{\g}.
\end{equation}
Applying $d^{\star}$ to the second equation above it follows easily that 
$$ L^{\star}_{d\theta}T^{\g}=0.
$$
provided that $n \neq 4$. The proof of the claim in (i) is obtained now by applying Lemma \eqref{extr} to the $2$-form $d\theta$.\\
(ii) By (i) the $1$-form $\theta$ is closed hence the curvature tensor $R^D$ belongs to 
$\mathcal{K}(\g)$ at any point of $M$, after also using (i) in Proposition \ref{Vvec}. If $\mathcal{K}(\g)=0$ then 
the connection $D$ must be flat. Otherwise $(\g,V)$ is again irreducible by Proposition \ref{grpvalg} and thus 
$\mathcal{K}(\g)$ is $1$-dimensional, generated (at some point of $M$) by a $\g$-invariant curvature tensor $R^{\g}$. The normalisation of $R^{\g}$ to have 
Ricci tensor equal to the metric ensures its $G$-invariance so by parallel transport we extend $R^{\g}$ over $M$ such that $DR^{\g}=0$. Hence 
\begin{equation} \label{curvlast}
R^D=\lambda R^{\g}
\end{equation}
for some smooth function $\lambda$ on $M$, which moreover turns out to be non identically zero. Indeed the vanishing of $\lambda$ would imply 
that of $R^D$ and hence the vanishing of $\g$ by using the Ambrose-Singer holonomy theorem, a contradiction with the irreducibility of $(\g,V)$. Now, the differential Bianchi identity for the connection $D$ reads (see \cite{besse}):
\begin{equation*}
\sigma_{X,Y,Z}(D_XR^D)(Y,Z)+R^D(Tor(X,Y),Z)=0
\end{equation*}
for all $X,Y,Z$ in $TM$, where $\sigma$ denotes the cyclic sum. Because the torsion tensor of $D$ is given by $Tor(X,Y)=\theta(X)Y-\theta(Y)X$ for all $X,Y$ in $TM$ we obtain 
$$ \sigma_{X,Y,Z}R^D(Tor(X,Y),Z)=2 \sigma_{X,Y,Z} \theta(X)R^D(Y,Z)
$$
whenever $X,Y,Z$ in $TM$. By also using \eqref{curvlast} and the fact that $DR^{\g}=0$ we arrive at 
$$ \sigma_{X,Y,Z} (d\lambda+2\lambda \theta )X \cdot R^{\g}(Y,Z)=0
$$
for all $X,Y,Z$ in $TM$, which is easily seen to imply that 
$d\lambda+2\lambda \theta=0$. Using Lemma \ref{lvglob}, it follows that $\lambda$ is nowhere zero and $\theta=-d ln \sqrt{\vert \lambda \vert}$. 
Let us consider now the Riemannian metric $\hat{g}=\varepsilon^2g $, where $\varepsilon=\sqrt{\vert \lambda \vert}$. 
We will show that $(M, \hat{g})$ is a locally symmetric space by examination of the conformal transformation rules for the relevant objects. First of all, using the 
Koszul formula for the Levi-Civita connections 
$\nabla$ and $\nabla^{\hat{g}}$ one easily arrives at 
\begin{equation} \label{confcomp}
\nabla^{\hat{g}}_XY=D_XY-\theta(X)Y
\end{equation}
for any smooth vector fields $X,Y$ on $M$, fact which leads after a straighforward calculation to $R^{\hat{g}}=R^D$. Therefore \eqref{confcomp} yields 
$$
(\nabla_X^{\hat{g}} R^{\hat{g}})(Y,Z)=(D_XR^D)(Y,Z)+2\theta(X)R^D(Y,Z)
$$ 
for all $X,Y,Z$ in $TM$. On the other hand side, since $DR^{\g}=0$ we have by using \eqref{curvlast}:
\begin{equation*}
\begin{split}
(D_XR^D)(Y,Z)=(X\lambda)R^{\g}(Y,Z)=2(X ln \varepsilon) R^D(Y,Z)=-2\theta(X)R^D(Y,Z)
\end{split}
\end{equation*}
for all $X,Y,Z$ in $TM$. We have showed that $\nabla^{\hat{g}} R^{\hat{g}}=0$, in other words $(M^n,\hat{g})$ is a locally symmetric space. \\
(iii) is an immediate consequence of the fact that $\theta$ is closed and of the behaviour of the holonomy algebra of $D$, after operating a conformal change in the 
metric, as the calculations in (ii) show. 
\end{proof}
An alternative partial proof of (i) in the Theorem above builds on the observation that $d\theta \wedge T^{\g}=0$, as it easily follows by differentiating the first equation in \eqref{strveceqs}. Indeed, it is enough 
to prove the following fact, well known for orthogonal representations on the Berger list.
\begin{pro} \label{multisymp}
Let $(G,V)$ be an irreducible orthogonal representation of some Lie group $G$, with Lie algebra denoted by $\g$. 
If $n=dim_{\mathbb{R}}V \geq 9$ and $T^{\g} \neq 0$ the map $L^{\g} : \Lambda^2V \to \Lambda^6V$ given by exterior multiplication with $T^{\g}$ is injective. 
\end{pro}
\begin{proof}
If $\alpha=\langle F \cdot, \cdot \rangle$ belongs to $Ker(L^{\g})$ a direct computation shows that 
\begin{equation*}
\begin{split}
0=L^{\star}_{\alpha}(\alpha \wedge T^{\g})=\vert F\vert^2 T^{\g}-2 \sum \limits_{i=1}^{n} Fe_i \wedge (Fe_i \lrcorner T^{\g})+2\alpha \wedge L^{\star}_{\alpha}T^{\g}
\end{split}
\end{equation*}
where $\{e_i, 1 \le i \le n \}$ is some orthonormal basis in $V$. Because our representation is irreducible and $T^{\g}$ is $G$-invariant we have $\langle X \lrcorner T^{\g}, Y \lrcorner T^{\g}\rangle=
\mu \langle X,Y\rangle  $ for all $X,Y$ in $V$ where $\mu=\frac{4}{n}\vert T^{\g}\vert^2$. Taking now above the scalar product with $T^{\g}$ yields 
$$0=\vert F\vert^2 \vert T^{\g} \vert^2-2\mu \vert F\vert^2+2\vert  L^{\star}_{\alpha}T^{\g} \vert^2=\vert F\vert^2 \vert T^{\g} \vert^2(1-\frac{8}{n})+2\vert  L^{\star}_{\alpha}T^{\g} \vert^2
$$
and the claim follows by using a positivity argument.
\end{proof}
The $8$-dimensional case is more involved and can be treated using structure results on 
isotropy algebras of $4$-forms in dimension $8$ from \cite{morgan, jlt2}. However the case of minimal representations is not covered, making the argument in Lemma \ref{symspace} necessary. 
When $n=4$, any Hermitian structure is preserved by a connection with vectorial torsion, determined by the (non-necessarily closed) Lee form of the structure. Therefore Theorem \ref{vec} cannot be extended over the $4$-dimensional case. 
\section{Some facts peculiar to dimension $8$}
Let $(\g,V)$ be a faithful orthogonal representation. We seek here to obtain a description of representations having the property that the map 
\begin{equation} \label{univG}
\varepsilon^{\perp} : \Lambda^3V \to \Lambda^1V \otimes \g^{\perp} \ \mbox{is surjective}.
\end{equation}
Under the additional assumption that $\g$ has the involution property, which includes the case when $\g$ is a simple Lie algebra it was shown 
in \cite{Fri2} that \eqref{univG} forces the equality $(\g,V)=(\mathfrak{spin}(7), \mathbb{R}^8)$.  We shall treat in this section, from a different perspective the 
general case. We observe first that: 
\begin{lema}
Any representation $(\g,V)$ satisfying \eqref{univG} must be irreducible.
\end{lema}
\begin{proof}
As observed before, if \eqref{univG} holds then
\begin{equation*}
(\Lambda^1V \otimes \g^{\perp}) \cap Ker(a)=0.
\end{equation*}
In particular, it follows that the map 
$$ a : S^2 \g^{\perp} \to \Lambda^4V
$$
is injective. Let us suppose now that we have an orthogonal splitting $V=V_1 \oplus V_2$ which is also $\g$-invariant.  Then $V_1 \wedge V_2 \subseteq \g^{\perp}$ and since 
$$ a((v_1 \wedge v_2) \otimes (v_1 \wedge v_2))=0
$$
for all $v_k$ in $V_k, k=1,2$ we obtain a contradiction. 
\end{proof}
We also need to establish the following general fact.
\begin{pro} \label{intertwin}
Let $(\g,V)$ be orthogonal and faithful such that $\Lambda^4V \cap (\Lambda^2V \otimes \g)=\{ 0 \} $. The map 
$$ I : (\mathcal{K}(\g))^{\perp} \cap Ker(Ric) \subseteq S^2(\g) \to S^2(\g^{\perp}) $$
given by 
$$ I(S)\alpha=L^{\star}_{\alpha}a(S)$$
for all $\alpha$ in $\g^{\perp}$ is well defined and injective.
\end{pro}
\begin{proof}Let us show first that $I$ is well defined. If $S$ is in $S^2(\g)$ a short computation using \eqref{formal} shows that 
\begin{equation*}
L^{\star}_{\alpha} a(S)=S \alpha-\{Ric(S), \alpha\}-\sum \limits_{x_k \in \g} [x_k, [Sx_k, \alpha]]
\end{equation*}
for all $\alpha$ in $\g$, where $\{x_k\}$ is some orthonormal basis in $\g$. If moreover $Ric(S)=0$, it follows that $L^{\star}_{\alpha} a(S)$ belongs to $\g$ for all $\alpha$ in $\g$ and one sees that 
$I(S)$ belongs to $S^2(\g^{\perp})$ by using an orthogonality argument. 

Let now $S$ in $S^2 \g \cap Ker(Ric)$ be such that $I(S)=0$. Then $L^{\star}_{\alpha}a(S)=0$ for all $\alpha$ in $\g^{\perp}$, and since 
$L^{\star}_{\alpha}a(S) $ belongs to $\g$ for all $\alpha$ in $\g$ it follows that $a(S)$ is in $\Lambda^4V \cap (\Lambda^2V \otimes \g)$, hence $a(S)=0$. But then 
$S$ must belong to $\mathcal{K}(\g)$ and so it vanishes. 

\end{proof}
We can now prove the following.
\begin{teo} \label{proofuniv}
Let $(\g,V)$ be an orthogonal and faithful representation such that \eqref{univG} holds. Then $(\g,V)=(\mathfrak{spin}(7), \mathbb{R}^{8})$.  
\end{teo}
\begin{proof}
Since our representation must be irreducible, we shall treat the various occurrences for $\mathcal{K}(\g)$ as Theorem \ref{clascurv} indicates. Also note that the 
dimension of the Lie algebra $\g$ is given by 
\begin{equation} \label{dimG}
\dim_{\mathbb{R}} \g^{\perp}=\frac{1}{6}(n-1)(n-2).
\end{equation}
since $\varepsilon^{\perp}$ is an isomorphism, in particular $\g$ is a proper subalgebra of $\mathfrak{so}(V)$.
If $\mathcal{K}(\g)=\{0\}$ or if it is $1$-dimensional, then we must have 
$\Lambda^4V \cap (\Lambda^2V \otimes \g)=\{0\}$ for $n \ge 4$, as it follows from 
Theorem \ref{class}, (ii). Now by Proposition \ref{intertwin}, the map $ I : S^2\g \cap Ker(Ric) \to S^2(\g^{\perp})$ is injective, hence 
$$dim_{\mathbb{R}} S ^2\g \le dim_{\mathbb{R}} S^2\g^{\perp}+\dim_{\mathbb{R}}S^2V \leq dim_{\mathbb{R}} S^2(\g^{\perp} \oplus V).$$
It follows that 
$ dim_{\mathbb{R}} \g \leq dim_{\mathbb{R}} \g^{\perp}+n$, which leads to $n^2 \leq 4+3n$ after using \eqref{dimG}. Thus $n \leq 4$, and the equality case can be excluded 
because then $\g$ would be $5$-dimensional and 
$$dim_{\mathbb{R}} S ^2\g=15  \nleq dim_{\mathbb{R}} S^2\g^{\perp}+\dim_{\mathbb{R}}S^2V=1+10=11.$$
It remains to treat the entries 
present in the Berger list \eqref{bergerlist}. It is an easy exercise to see that all entries but the 6-th can be excluded for dimensional reasons.
Hence only the case $(\g,V)=(\mathfrak{spin}(7), \mathbb{R}^{8})$ is eligible, when moreover \eqref{univG} is known to hold \cite{Iv,Fri2}, and the proof is complete.  
\end{proof}
\begin{rema} Another interesting $\g$-module is $(\Lambda^1V \otimes \g) \cap Ker a \subseteq \Lambda^1V \otimes \g$. It appears naturally in the study of holonomy 
groups of Lorentzian metrics in connection with the so-called weak Berger algebras \cite{L1}. It seems not unlikely that techniques similar to those in this paper 
could lead to a characterisation of the instances when $(\Lambda^1V \otimes \g) \cap Ker(a)$ vanishes.
\end{rema}
$\\$
{\bf{Acknowledgements}}: This research has been partly supported by the Royal Society of New Zealand, Marsden Grant no. 06-UOA-029 and the 
University of Hamburg where the final part of this work has been 
carried out. It is a pleasure to thank R.Gover, Th.Leistner, P.Nurowski and especially  V.Cort\'es for useful discussions.

\end{document}